\documentclass[10pt,a4paper,reqno,twoside]{amsart}
\usepackage{amsmath}
\usepackage{amssymb}
\usepackage{amsthm}
\usepackage{mathrsfs}
\usepackage{color}
\usepackage{amsmath}
\usepackage{amsfonts}
\usepackage{amsthm}

 \newtheorem{definition}{Definition}[section]
  \newtheorem{example}{\sc Example}[section]
 \newtheorem{notation}{Notation}
 \newtheorem{hypothesis}{Hypothesis}
 \newtheorem{remark}{Remark}

\newcommand{\R}{\mathbbm R}

\newcommand{\bec}{\begin{cases}}
\newcommand{\eec}{\end{cases}}

\newcommand{\pd}{{\mathrm{pd}}}
\newcommand{\hyp}{{\mathrm{hyp}}}

\newcommand{\abs}[1]{\left\vert#1\right\vert}

\title[Wave models with time-dependent mass and speed of propagation.]{A classification for wave models with time-dependent mass and speed of propagation}
\author[M. R. Ebert\and W. N. Nascimento]{Marcelo Rempel Ebert \and Wanderley Nunes do Nascimento }

\address{Marcelo Rempel Ebert\\
	Department of Computing and Mathematics (FFCLRP),
	University of S\~ao Paulo (USP),
	Av. Bandeirantes, 3900 - CEP 14040-901 - Ribeir\~ao Preto - SP -Brazil}
\email{ebert@ffclrp.usp.br}

\address{Wanderley Nunes do Nascimento\\
	Department of Mathematics,	
	Institute of Mathematics, Statistics, and Computer Science (IMECC),	
	State University of Campinas (Unicamp),
Rua S\'ergio Buarque de Holanda, 651 - CEP
13083-859 - Campinas - SP - Brazil}
\email{wnunesmg@yahoo.com.br}

\newtheorem{thm}{Theorem}[section]
\newtheorem{lem}[thm]{Lemma}
\newtheorem{prop}[thm]{Proposition}

\theoremstyle{definition}

\theoremstyle{remark}
\newtheorem{rem}{Remark}[section]

\newtheorem{exam}[rem]{Example}
\numberwithin{equation}{section}

\def\d{\mathrm d}

\DeclareMathOperator{\real}{Re}
\DeclareMathOperator{\diag}{diag}

\def\R{\mathbb R}

\def\d{\mathrm d}
\let\Im\relax
\DeclareMathOperator{\Im}{Im}
\begin{document}

\begin{abstract}
In this paper, we study  the long time behavior of  energy solutions for a class of wave equation with time-dependent mass and speed of pro\-pagation.
We introduce a classification of the potential term, which clarifies whether the solution behaves like the solution
to the wave equation or Klein-Gordon equation.  Moreover,   $L^q-L^2, q\in [1, 2]$ estimates   for  scale-invariant models are derived and  applied to obtain
 global in time small data energy solutions  for the semilinear  Klein-Gordon
equation in anti de Sitter spacetime.
\end{abstract}
%%%%%%%%%%%%%%%%%%%%%%%%%%%%%%%%%%%%%%%%%%%%%%%%%%%%%%%%%%%%%%%%%%%%%
\thanks{The first author has been partially
supported by FAPESP Grant 2015/16038-2. The second author is supported by FAPESP  Grant 2015/23253-7.}
%\thanks{The first author thanks Michael Reissig and Marcelo R. Ebert for the values discussion about the subject. }
\maketitle
\section{Introduction  }

Let us consider the Cauchy problem for the wave equation with time-dependent mass and speed of propagation
\begin{equation}\label{wavegeneral}
\begin{cases}
u_{tt} - a(t)^2 \Delta u + m(t)^2u=0,
  \quad (t,x)\in (0,\infty)\times \R^n,
\\
(u(0,x),u_t(0,x))=(u_0(x),u_1(x)),
  \quad x\in \R^n.
\end{cases}
\end{equation}
The  Klein-Gordon type energy for the solution to \eqref{wavegeneral} is given by
\begin{equation}\label{E}
  E_{a,m}(t)\doteq\frac12\big( \|u_t(t,\cdot)\|_{L^2}^2 +  a(t)^2\|\nabla u(t,\cdot)\|_{L^2}^2 +  m(t)^2\| u(t,\cdot)\|_{L^2}^2\big).
\end{equation}
One can observe many different effects for the behavior of~$E_{a,m}(t)$ as~$t\to\infty$ according to the properties of the speed of propagation~$a(t)$ and the coefficient $m(t)$ in the mass term.

We first discuss  properties of the energy in the case $m(t)\equiv 0$ in \eqref{wavegeneral}.
If $0<a_0\leq a(t)\leq a_1$ for any~$t\geq0$, then the energy $E_{a,0}(t)$ is equivalent to
\[ E_1(t)=\frac12\big( \|u_t(t,\cdot)\|^2_{L^2}+ \|\nabla u(t,\cdot)\|^2_{L^2}\big). \]
Although $E_{1,0}(t)$ is a conserved quantity, oscillations of~$a(t)$ may have a very deteriorating influence on the energy behavior for the solution to~\eqref{wavegeneral} (see~\cite{C} and \cite{R-Y}). On the other hand, if~$a\in\mathcal{C}^2$ and
\[
 |a^{(k)}(t)|\leq C_k(1+t)^{-k} \quad \text{for $k=1,2$}
 \]
(so only very slow oscillations are allowed), then the so-called \emph{generalized energy conservation} property holds~\cite{ReiSmith}. This means, there exist positive constants $C_0$ and $C_1$ such that the inequalities
\begin{equation}\label{eq:GEC}
C_0E_{a,0}(0)\leq E_{a,0}(t)\leq C_1E_{a,0}(0)
\end{equation}
are valid for all $t \in (0,\infty)$, where the constants are independent of the data.\\
If  ~$a(t)\geq a_0>0$ is an increasing function  satisfying a suitable control on the oscillations,
then one can prove the estimates~\cite{ReisBui}
\begin{equation}\label{eq:HWGEC}
E_{a,0}(t)\leq C a(t)\,\big(E_{a,0}(0)+\|u_0\|_{L^2}^2\big).
\end{equation}
We remark that an essential difference between~\eqref{eq:GEC} and \eqref{eq:HWGEC}  is that on the right-hand side of~\eqref{eq:HWGEC} it appears the~$H^1$ norm of~$u_0$, not only the $L^2$-norm of its gradient as in~\eqref{eq:GEC}.

In the case $a(t)\equiv 1$,  $E_{1,m}(t)$ is a conserved quantity for the classical Klein-Gordon equation, whereas it is known that the behavior of the potential energy $\| u(t,\cdot)\|_{L^2}$ change accordingly $\lim_{t \to \infty} t m(t)=\infty$ or $\lim_{t \to \infty} t m(t)=0$. To explain this effect, let us consider the energy
\[E_{p}(u)(t) \doteq\; \frac12
\Big(\|u_t(t,\cdot)\|_{L^2}^2+\|\nabla_x u(t,\cdot)\|_{L^2}^2 + p(t)^2\|u(t,\cdot)\|_{L^2}^2\Big)\]
In the PhD thesis
\cite{B}, the author studied decreasing coefficients $m=m(t)$  which satisfy among other things $\lim_{t \to \infty} t m(t)=\infty$. In this case the potentials   are called  \textit{effective}, i.e., the decays of the solution and its derivatives are related with the decays of the classical Klein-Gordon equation measured in the $L^q$ norm. Under some additional condition on $m$, the following  energy estimate  was derived
\begin{equation}\label{energyest}
E_{p}(u)(t) \leq C E_{p}(u)(0),
\end{equation}
 with ${p(t)^2}=m(t)$. For decreasing $m$ estimate \eqref{energyest} is better than an estimate like $E_{1,m}(t)\leq C E_{1,m}(0)$.
 In \cite{B1},  the authors  derived  the energy estimate \eqref{energyest} for scale invariant models $m(t)=\frac{\mu}{1+t}, \mu>0$, but now the constant $\mu$ has  influence
 in the function $p(t)$.

  In \cite{EKNR, Nthesis} the authors explained qualitative properties of solution to \eqref{wavegeneral}
  in the case $\lim_{t \to \infty} t m(t)=0$. If $(1+t)m(t)^2 \in L^1(\R^+)$, it was proved a scattering result to
  free wave equation, whereas the potentials are called \textit{non-effective}, i.e., the decays of the solution and its derivatives are related with the decays of the free wave  equation measured in the $L^q$ norm. The energy estimate \eqref{energyest}, with $p(t)=(1+t)^{-1}\psi(t)$  and $\psi$ an increasing function, was obtained in the    case of  $(1+t)m(t)^2 \notin L^1(\R^+)$.

  In \cite{GalY} the authors derived an explicit solution  to  the
 Cauchy problem for the  well-known  Klein-Gordon
equation in anti de Sitter spacetime
\begin{equation}\label{deSitter}
u_{tt} - e^{2t} \Delta u  + m^2 u
=0,\,\,\,u(0,x)=u_0(x),\,\,\,u_t(0,x)=u_1(x).
\end{equation}
By using the obtained representation of solutions, in \cite{Gal} it was proved some  $L^q-L^{q'}$ estimates, with $q\in (1, 2]$ and $\frac1{q}+\frac1{q'}=1$,
 but exception  for the case
of space dimension $n=1$, we have some loss of regularity with respect to the initial data. Due to the lack of $L^q-L^r$ estimates without loss of regularity,
 with $1\leq q\leq r\leq \infty$,  it is a challenging problem to derive the  critical exponents for global (in time) small data
energy solutions  to the Cauchy problem for the semilinear Klein-Gordon equation in anti de Sitter spacetime
 \begin{equation}\label{CauchyProbN3}
 u_{tt} -  e^{2t}\Delta u  + m^2 u
=|u|^p,\,\,\,u(0,x)=u_0(x),\,\,\,u_t(0,x)=u_1(x).
\end{equation}

 The same difficulties took place in the treatment of the  classical Klein-Gordon semilinear equation  with power nonlinearity $|u|^p$,
	being that only in the late nineties it was shown that, for space dimension $n \leq 3$,
	the critical exponent is the well known Fujita index $p_{Fuj} \doteq 1 + \frac{2}{n}$, i.e.,
	global (in time) existence of small data energy solutions holds
	for   $p> p_{Fuj} $ (see \cite{LS}),  whereas blow up results are established in \cite{KT}
	for $1 < p\leq p_{Fuj}$.

In this paper, our main goal is to  introduce a classification (see Definition \ref{Def5.1}) for the potentials in \eqref{wavegeneral} in terms of the time-dependent speed of propagation $a(t)$, consistently extended from the case of constant speed of propagation $a\equiv 1$.
	In the case of \textit{effective} and \textit{non-effective} potentials we derive sharp energy estimates and we show optimality of the results
	by the aid of scale-invariant models.
	Moreover, we  explain some gaps from the thesis \cite{B} and  we derive $L^q-L^2$ estimates, with $q\in [1,2]$,  for the model \eqref{deSitter}.
	As an application to our derived linear estimates, we proved global existence (in time) of small data energy solutions for the Cauchy problem
	\eqref{CauchyProbN3}.
% In Theorem \ref{nonlinear}, at least for dimension $n=1,2$, we have some improvement for $p$, this explain us that by replacing the speed of propagation of the classical Klein-Gordon equation  by  an exponential speed of propagation produce an additional dissipative
% effect in the equation.

For the ease of reading, we summarize the scheme of the paper:
\begin{itemize}
	
	\item in Section \ref{sectmain}, we propose a classification for the potential term and we state our main  results;

	\item in  Section \ref{diagprocedure}, we describe a diagonalization procedure to be used in sections  \ref{secteffect} and \ref{noneffect};
	
	\item in  Section \ref{secteffect}, we prove the result for effective potential;
	
	\item in Section \ref{noneffect}, we prove the result for non-effective potential;
	
	\item in Section \ref{scattering}, we prove the Scattering result;
	
	\item in Section \ref{scaleinvariant}, we discuss some scale invariant models and we prove Theorem \ref{nonlinear};

	\item Section \ref{open} completes the paper with concluding remarks and  open problems;
	
\end{itemize}

\section{Main Results}\label{sectmain}
In this paper, we use the following notations.
\begin{notation}
Let $f,g : \Omega \subset \mathbb{R}^n \to\R$ be two strictly positive functions. We use the notation $f \approx g$ if there exist two constants $C_1,C_2>0$ such that $C_1 g(y)\leq f(y) \leq C_2 g(y)$ for all $y\in \Omega$. If the inequality is one-sided, namely, if $f(y)\leq Cg(y)$ (resp.~$f(y)\geq Cg(y)$) for all~$y\in \Omega$, then we write~$f\lesssim g$ (resp.~$f\gtrsim g$).
\end{notation}

To state our results and assumptions on the coefficients of the equation in~\eqref{wavegeneral} we introduce some auxiliary functions.
\begin{definition}
Let~$a\in\mathcal{C}^2[0, \infty)$ be a strictly positive function. We define
\[
A(t)  \doteq 1+ \int_0^t a(\tau)d\tau, \quad
\eta(t)  \doteq \frac{a(t)}{A(t)}.\]
%Let~$m\in\mathcal{C}^2[0, infty)$ be a strictly positive function. We define
%\[
%M(t)  \doteq 1+ \int_0^t m(\tau)d\tau, \quad
%\rho(t)  \doteq \frac{m(t)}{M(t)}.\]
%We also define
%\[ \gamma(t) \doteq \max \{\eta(t), \rho(t)\}.\]
\end{definition}
To study the interaction between $a$ and $m$ we assume  the following conditions:
\begin{hypothesis}\label{Hyp1}
We assume that~$a\in\mathcal{C}^2[0, \infty)$, $a(t)>0$, with $a\not\in L^1$ and $a(0)=1$,   together with the estimates
\begin{equation}\label{shape}
\frac{|a^{(k)}(t)|}{a(t)} \lesssim  \eta(t)^k\, \ \text{for $k=1,2$}.
\end{equation}
\end{hypothesis}
\begin{hypothesis}\label{Hyp2}
We assume that~$m\in\mathcal{C}^2[0, \infty)$  has the form
\[m(t)=\mu(t)\eta(t)> 0,\]
with $m(0)=1$,  may have an oscillating behavior.  For this reason we suppose
\begin{equation}\label{shapeOsc}
|\mu^{(k)}(t)| \lesssim \, \mu(t) \eta(t)^k\, \ \text{for $k=1,2$}.
\end{equation}
\end{hypothesis}
\begin{definition} \label{Def5.1}
We propose the following classification of potential terms:
\begin{enumerate}
\item The potential term $m(t)^2 u$ generates \emph{scattering} to the corresponding wave model if
\begin{align} \label{ScatCond}
	\frac{A(t)}{a(t)}m(t)^2 \in L^1(\R^n).
	\end{align}

\item  The potential term $m(t)^2 u$ represents a \emph{non-effective potential}
if
\begin{equation}\label{noneffective}
\lim_{t \to \infty}\mu(t)=0,
\end{equation}
but $	\frac{A(t)}{a(t)}m(t)^2 \notin L^1[0,\infty).$
	
\item The potential term $m(t)^2 u$ generates  an \emph{effective potential} if
\begin{equation}\label{effectivecondition}
\lim_{t \to \infty} \mu(t) =\infty.
\end{equation}
\end{enumerate}
\end{definition}
But there exists a grey zone. The models of the
grey zone are in the boundary between effective and non-effective potentials, they  can be described by the models
%scale invariant models
\[
u_{tt} - a(t)^2\triangle u +  \mu^2 \frac{a(t)^2}{A(t)^2}u =0 \,,\]
where $ \mu\geq 0$ is a constant. For these models,  the size of the constant
$\mu$ may have some  influence in the long time behaviour,  namely, for small $\mu$ we are in the case of non-effective potentials and we are able to include it in Theorem \ref{NoneffTeo}, whereas for large $\mu$ we are in the case of
effective potential, but for simplicity of the proof we did not included it in Theorem \ref{effectiveth}. For this reason,
we discuss some scale invariant models in sections \ref{scaleinvariant}.

\subsection{Effective potential}

To  state our result in the case of effective potential we define
\begin{eqnarray} \label{decayrate}
\gamma(t)\doteq \max \{a(t), m(t)\}
\end{eqnarray}
and the following energy
\begin{eqnarray} \label{Energy}
E(u) (t)\doteq \frac{1}{2} \left(\|u_t(t,\cdot) \|_{L^2}^2 + a(t)^2\|\triangledown_x u(t,\cdot) \|_{L^2}^2 + m(t)\gamma(t)\|u(t,\cdot) \|_{L^2}^2 \right).
\end{eqnarray}
\begin{thm}\label{effectiveth}
 Let $(u_0,u_1) \in H^1 \times L^2$ and $u$ be an energy solution of the Cauchy problem \eqref{wavegeneral}. We assume condition \eqref{effectivecondition}, Hypotheses \ref{Hyp1} and \ref{Hyp2}.  In addition, if $\frac{\eta}{\mu}\in L^1[0, \infty)$, then we have the following estimate for the energy
	\begin{equation}\label{effectiveRate}
	E(u)(t) \lesssim \gamma(t) E(u)(0), \qquad \forall t\geq 0,
	\end{equation}
where $\gamma$ is given by \eqref{decayrate}.
\end{thm}
\begin{remark} If $\frac{a}{m}$ is bounded, then
\[ \frac{\eta}{\mu}=\frac{a}{m}\frac{a}{A^2}\in L^1.\]
 If  $\frac{\eta}{\mu}\notin L^1[0, \infty)$, we are near to models in the so-called  grey zone.
 \end{remark}
\begin{remark}
 From  Theorem \ref{effectiveth} we conclude that for increasing $m$  the potential energy decay, i.e.,
\[\|u(t,\cdot) \|_{L^2}^2\lesssim \frac1{ m(t)} E(u)(0), \qquad \forall t\geq 0.\]
This estimate  is better than the conjecture done in \cite{B}.
\end{remark}
\begin{remark} If $\max \{a(t), m(t)\}=m(t)$, under additional regularity on the initial data and using the a priori estimates for $\|u(t,\cdot) \|_{L^2}$, one may derive a better decay for the elastic energy, namely,
\[a(t)^2\|\triangledown_x u(t,\cdot) \|_{L^2}^2  \lesssim a(t) ( \|u_0 \|_{H^2}^2+ \|u_1 \|_{H^1}^2).\]
\end{remark}
\begin{exam}
Let $a(t)=m(t)$, that is, take $\mu(t)=A(t)$. Thanks to $a\not\in L^1$,  it is clear that $\frac{\eta(t)}{\mu(t)}=\frac{a(t)}{A(t)^2}\in L^1[0, \infty)$ and the conclusion of
Theorem \ref{effectiveth} holds with $\gamma(t)=a(t)$. In particular, for increasing $a$ the potential energy decay
\[\|u(t,\cdot) \|_{L^2}^2\lesssim \frac1{ m(t)} E(u)(0), \qquad \forall t\geq 0.\]
In \cite{ReiYad} the authors derived $L^p-L^q$ estimates for the elastic and kinetic energy for this model.
\end{exam}
\begin{exam}
Let $a(t)=(1+t)^{\ell}$, with $\ell>-1$, and $m(t)=\mu(1+t)^{\epsilon -1}$, with $\epsilon>0$. Then $\mu(t)\sim(1+t)^{\epsilon}$ and the statement of Theorem \ref{effectiveth} holds with
$\gamma(t)=(1+t)^{\ell}$ for $\epsilon\leq \ell+1$ and $\gamma(t)=(1+t)^{\epsilon -1}$ for $\epsilon> \ell+1$. (see Example \ref{polinomial} for the limit case $\epsilon=0$)
\end{exam}
\begin{exam}
Let $a(t)=e^{t}$ and $m(t)=\mu(1+t)^{\epsilon }$, with $\epsilon>1$.  Then the statement of Theorem \ref{effectiveth} holds with
$\gamma(t)=e^{t}$.
\end{exam}

\subsection{Non-effective potential}

\

To  state our result in the case of non-effective potential we assume   the following hypotheses:
\begin{hypothesis}\label{Hyp3}
There exists a positive non-decreasing function $\psi\in C^2(\R^+)$ with $\psi(0)=1$ such that
$\frac{1}{\eta(t)\psi(t)^2}$ is
increasing for large $t$ and
\begin{equation}
\label{neweq32}  \frac{|\psi'(t)|}{\psi(t)} <c \eta(t), \qquad c\in (0,1).
\end{equation}
 Besides \eqref{neweq32} the following relation between $m$, $\eta$ and $\psi$  must be satisfied$:$
\begin{equation}\label{neweq33}
\psi(t)^2\eta(t)\int_{0}^t\psi(\tau)^{-2}d\tau + \int_0^{\infty}\frac1{\eta(\tau)}\Big|\frac{\psi''(\tau)}{\psi(\tau)}+m(\tau)^2\Big|d\tau \lesssim 1.
\end{equation}
\end{hypothesis}
Now  we define
\begin{eqnarray}
\label{EnergNonEff}
E_{p}(u)(t) = \frac12 \left( \| u_t(t,\cdot) \|_{L^2}^2+ a(t)^{2} \| \nabla u(t,\cdot) \|_{L^2}^2 + p(t)^2\| u(t,\cdot) \|_{L^2}^2 \right),
\end{eqnarray}
with
\[p(t) \doteq  \eta(t) \psi(t) \sqrt{ q(t)}, \qquad q(t)\doteq \max \{a(t), \psi(t)^{-2} \}.\]
%and
%\begin{eqnarray}
%\label{DefPsiNon}
%\psi(t) \doteq \exp\left( \int_0^t  \frac{\mu(\tau)^2}{1+\tau} d\tau\right).
%\end{eqnarray}
 Then we have the following energy estimate:
%%%%%%%%%%%%%%%%%%%%%%%%%%%%%%%%%%%%%%%%%%%%%
\begin{thm}\label{NoneffTeo}
	Let $(u_0,u_1) \in H^1 \times L^2$ and $u$ be the solution of the Cauchy problem \eqref{wavegeneral}. We assume Hypotheses  \ref{Hyp1} to \ref{Hyp3}.
 Then we have the estimate
	\begin{eqnarray}\label{energyinequality}
	E_{p}(u)(t) \lesssim q(t) E_p(u)(0), \qquad \forall t \geq 0.
	\end{eqnarray}
\end{thm}
\begin{remark}
As one may verified in the next examples, the condition \eqref{noneffective} takes place in order to guarantee  the existence of the function $\psi$ in Hypothesis  \ref{Hyp3}.
\end{remark}
%%%%%%%%%%%%%%%%%%%%%%%%%%%%%%%%%%%%%%%%%%%%%%%%%%%%%%%
\begin{example}(Scale invariant model)\\
Consider the  Cauchy problem
\[
u_{tt} - (1+t)^{2\ell} \Delta u  + \frac{\mu^2}{(1+t)^2} u
=0,\,\,\,u(0,x)=u_0(x),\,\,\,u_t(0,x)=u_1(x),
\]
where $\ell>-1$ and $0<\mu^2 < \frac1{4}$.  Let us take the function $\psi$
   from Hypothesis \ref{Hyp3} as
  \[\psi(t)=(1+t)^{\sigma}, \qquad {\text with } \qquad 2\sigma=1-\sqrt{1-4\mu^2}.\]
  It is clear that $\frac{\psi''(\tau)}{\psi(\tau)}+m(\tau)^2=0$ and all the conditions from Hypotheses \ref{Hyp1} to  \ref{Hyp3} are satisfied.
If $(u_0,u_1) \in H^1 \times L^2,$ then Theorem \ref{NoneffTeo} implies the following estimates
\[
\|u(t,\cdot) \|_{L^2}^2 \lesssim
(1+t)^{1+\sqrt{1-4\mu^2}} \]
and
\begin{eqnarray*} \label{eq: p(t)}
\|u_t(t,\cdot) \|_{L^2}^2 +(1+t)^{2\ell}\|\nabla u(t,\cdot) \|_{L^2}^2 \lesssim \left\{ \begin{array}{cr}
(1+t)^{\ell}, & \ell+1> \sqrt{1-4\mu^2},\\
(1+t)^{-1+\sqrt{1-4\mu^2}},& \ell+1\leq \sqrt{1-4\mu^2}. \end{array} \right.\end{eqnarray*}
\end{example}
\begin{rem}
The previous example can be treated  in a different way, including also large parameter $\mu\geq \frac1{4}$, see Example \ref{polinomial}.
\end{rem}
\begin{example}\label{oldth}
	Let $a(t)=(1+t)^{\ell}$ with $\ell>-1$. If $m(t)=\frac{\mu(t)}{e+t}$, with $\mu(t)$ satisfying
conditions \eqref{shapeOsc}, \eqref{noneffective} and
$
\frac{\mu(t)^4}{e+t} \in L^1,$  then the
conclusion of Theorem \ref{NoneffTeo} holds with $\psi(t)$ given by
\[\psi(t) \doteq \exp\left( \int_0^t  \frac{\mu(\tau)^2}{e+\tau} d\tau\right).\]
Indeed,
\[ \frac{\psi''(\tau)}{\psi(\tau)}+m(\tau)^2= \frac{2\mu(t)\mu'(t)}{e+t}+ \frac{\mu(t)^4}{(e+t)^2}  \]
Moreover, condition \eqref{noneffective} implies that
\begin{equation}\label{newWirthLemma}
\int_{0}^t\psi(s)^{-2}ds\approx\frac{t}{\psi(t)^2},
\end{equation}
and $\frac{t}{\psi(t)^2}$ is increasing for large $t$.
Indeed, integration by parts yields
\[ \int_{0}^t\psi(s)^{-2}ds = \frac{t}{\psi(t)^2} + 2\int_0^t \frac{s}{\psi(s)^2}  \frac{\mu(s)^2}{e+s}ds. \]
On the one hand the right side is greater than $t\psi(t)^{-2}$. On the other hand,\\
 $2t(e+t)^{-1}\mu(t)^2 < \epsilon$  for large $t$ and any $\epsilon>0$. Then,
$$ \int_{0}^t\psi(s)^{-2}ds \leq \frac{1}{1-\epsilon}\left( \frac{t}{\psi(t)^2} \right) \lesssim \frac{t}{\psi(t)^2}.  $$
Monotonicity is a consequence of
\[  \frac{d}{dt}\frac{t}{\psi(t)^2} =\left( 1-\frac{2t\mu(t)^2}{(1+t)}  \right) \frac{1}{\psi(t)^2},   \]
which is positive for large $t$. Therefore, all the conditions from Hypotheses \ref{Hyp1} to  \ref{Hyp3} are satisfied.\\
For instance, if   $\mu(t)^2 = \frac{\mu^2}{\ln(e+t)\cdots
		\ln^{[k]}(e^{[k]}+t)} $  with $e^{[j+1]}=e^{e^{[j]}}$ and
	$\ln^{[j+1]}(t)=\ln(\ln^{[j]}(t))$, $j=1,2, \cdots$, then
	$\psi(t)\approx (\ln^{[k]}(e^{[k]}+t))^{\mu^2}$.\\
If   $\mu(t)=\frac{\mu}{\ln(e+t))^{\gamma}}$, with $\mu>0$ and  $\gamma>\frac1{4}$, then  $\psi(t)$ given by
\[
\psi(t)=\begin{cases}
\exp (\mu^2(\ln (e+t))^{1-2\gamma}), \gamma \in (\frac1{4}, \frac1{2}) \,,\\
\ln (e+t))^{\mu^2},  \gamma=\frac1{2}\,,\\
1, \gamma>\frac1{2}.
\end{cases}
\]
\end{example}
\begin{rem}\label{remcatalanNumber}
In Example \ref{oldth}, if $\frac{\mu(t)^4}{1+t} \notin L^1$, one  may still derive a  result as in Theorem \ref{NoneffTeo} (see \cite{EKNR}) by replacing $\psi$ by
\begin{equation}
\label{modeleq4} \psi(t)=\exp\Big(\sum_{k=1}^{N}\gamma_k\int_0^t\frac{\mu(\tau)^{2k}}{(1+\tau)}d\tau\Big),
\end{equation}
for some integer $N>1$. The sequence $\{\gamma_k\}_k$ in \eqref{modeleq4} is relate to the well-known   Catalan numbers which
 can be found for instance in  \cite{K}.
\end{rem}

\subsection{Scattering to wave equation}

\

In this section we will impose conditions for $a(t)$ and  $m(t)$ such that the
solutions $u=u(t,x)$ of Cauchy problem  \eqref{wavegeneral}
behave asymptotically equal to the solution of corresponding wave equation with strictly increasing speed of propagation
\begin{eqnarray}\label{CauchyProbSpeeWave}
v_{tt} - a(t)^2 \Delta v
=0,\,\,\,v(0,x)=v_0(x),\,\,\,v_t(0,x)=v_1(x),
\end{eqnarray}
with some suitable Cauchy data $(v_0,v_1)$.

Let us define the function space
$$  E = {L^2(\R^n) \times L^2(\R^n)}.$$
Before stating the result we define for any~$\epsilon>0$ the
following closed subset of~$E$:
\[ F_\epsilon := \left\{ U_0\in E : \ U_0(\xi)=0 \ \text{for any~$|\xi|\leq\epsilon$} \right\}. \]
We remark that~$\mathcal{L}=\cup_{\epsilon>0}F_\epsilon$ is a dense
subset of~$E$.

In addition to Hypothesis \ref{Hyp1}, we assume  $a'(t)>0$  for all $t\in [0,\infty)$,   together with the estimate
\begin{equation}\label{additionalscat}
\sqrt{a(t)}
\int_{0}^t \sqrt{a(\sigma)}d\sigma \lesssim A(t).
\end{equation}

%%%%%%%%%%%%%%%%%%%%%%%%%%%%%%%%%%%%%%%%%%%%%%%%%%%%%%%%
%%%%%%%%%%%%%%%%%%%%%%%%%%%%%%%%%%%%%%%%%%%%%%%%%%%%%%%
\begin{thm}	\label{scat}
	Let us  assume Hypothesis \ref{Hyp1}, conditions  \eqref{ScatCond},  \eqref{additionalscat}  and that    $a'(t)> 0$  for all $t\in [0,\infty)$.
		Then, for any initial data $(u_0,u_1)\in H^1\times L^2$, there exists a linear, bounded operator~$W_+(D):E\to E$
	such that if the
	initial data of the Cauchy problems \eqref{wavegeneral}	and \eqref{CauchyProbSpeeWave} are related by $(a(0)\nabla v_0,v_1)=W_+(D) (\langle D(0)\rangle u_0,u_1)$, it follows that  the asymptotic equivalence of solution  holds
	\begin{equation}
	\label{eq:modscatt}
	\lim_{t\to\infty} \frac{1}{\sqrt{a(t)}}  \Big\|(a(t)\nabla v(t,\cdot),v_t(t,\cdot))-\Big(\langle D(t)
	\rangle u(t,\cdot),u_t(t,\cdot)\Big)\Big\|_{E} = 0\,,
	\end{equation}
	where $\langle D(t)\rangle$ denotes the pseudodifferential operator having the symbol
	\begin{equation}\label{hfunction}
	h(t,\xi)=\left(|\xi|^2a(t)^2+N^2 \eta(t)^2\right)^{\frac{1}{2}}.
	\end{equation}
	Moreover, on the dense subset $\mathcal{L}$ we can
	state the decay rate as
	\begin{equation}
	\label{eq:modscatt2}
	\frac{1}{\sqrt{a(t)}} \Big\|(a(t)\nabla v(t,\cdot),v_t(t,\cdot))-\Big(\langle D(t)
	\rangle u(t,\cdot),u_t(t,\cdot)\Big)\Big\|_{\mathcal{L}} \lesssim \int_t^\infty
	\frac{A(\tau)}{a(\tau)}m^2(\tau)d\tau
	\end{equation}
	as $t$ goes to infinity.
\end{thm}
\begin{rem}
Condition \eqref{additionalscat} holds for a large class of examples like $a(t)=(1+t)^{\ell}, \ell \geq 0$, $a(t)=e^t$ and $a(t)=e^{e^t}e^t$.
 More in general, it is true under the assumption   $a'(t)\leq C_1 a(t)\eta(t)$ for all $0<C_1<2$ in   \eqref{shape}.
Indeed, let us consider the function
\[F(t)=\int_{0}^t \sqrt{a(\sigma)}d\sigma-C \frac{A(t)}{\sqrt{a(t)}}.\]
Then
\[ F'(t)=\sqrt{a(t)} - C\sqrt{a(t)} + \frac{C A(t)}{2\sqrt{a(t)}} \frac{a'(t)}{a(t)}\leq  \left(1-C + \frac{CC_1}{2}\right)\sqrt{a(t)}<0, \]
for $C>0$ sufficiently large.
\end{rem}
\subsection{Semilinear  Klein-Gordon
equation in anti de Sitter spacetime}

\

Let us consider the  Cauchy problem for the semilinear  Klein-Gordon
equation in anti de Sitter spacetime
\begin{equation}\label{CauchyProbN2}
u_{tt} - e^{2t} \Delta u  + m^2 u
=|u|^p,\,\,\,u(0,x)=u_0(x),\,\,\,u_t(0,x)=u_1(x),
\end{equation}
with $m>0$ and $p>1$.
Let us define the function spaces
$$ \mathcal{D}^\kappa_q(\mathbb{R}^n) = (H^{\kappa + 1} \cap L^q) \times (H^{\kappa} \cap L^q)$$
with $q \in [1,2)$ and the norm $\| (u,v) \|_{\mathcal{D}^\kappa_q}^2 = \| u \|_{L^q}^2+\| u \|_{H^{\kappa+1 }}^2+\| v \|_{L^q}^2+\| v \|_{H^{\kappa}}^2.$
In the following, we denote $\mathcal{D}^0_q=\mathcal{D}_q$.
\begin{thm}\label{nonlinear} Let $n\leq 4$,  $m>0$ and
\[  2\leq p\leq \frac{n}{[n-2]_+}.\]
Then there exists a constant $\epsilon>0$ such that for all $(u_0,u_1)\in \mathcal{D}_1(\mathbb{R}^n)$ with
\[||(u_0,u_1)||_{\mathcal{D}_1(\mathbb{R}^n)}\leq \epsilon\]
there exists a uniquely determined energy solution $u\in C([0, \infty), H^1(\mathbb{R}^n))\cap C^1([0, \infty), L^2(\mathbb{R}^n)$ to  \eqref{CauchyProbN2}.
Moreover,  the solution satisfies the following estimates
\[
\| u_t(t,\cdot)\|_{L^2}+ e^t\|\nabla_x u(t,\cdot)  \|_{L^2} \lesssim e^{t/2}   \big( \| u_0 \|_{H^1}+ \| u_1\|_{L^2} \big), \qquad  \forall t\geq0,
\]
\[
\|u(t,\cdot)\|_{L^2} \leq C\,e^{-t/2}d(t)\,\|(u_0,u_1)||_{\mathcal{D}_1(\mathbb{R}^n)}, \qquad  \forall t\geq0.
\]
where
\begin{align*}
		 d(t)= \left\{ \begin{array}{cr}
		1  &\mbox{\,\,\, for \,\,\,}  n\geq 2,\\
		t^{\frac{1}{2}}   &\mbox{\,\,\, for \,\,\,} n=1.\end{array} \right.
		\end{align*}
\end{thm}
\begin{rem} \label{embedding}
If $n=1$, by using the embedding of $H^1(\R)$ into $L^{\infty}(\R)$  and interpolation results, we  no longer need to use
Gagliardo-Nirenberg inequality in the proof of Theorem \ref{nonlinear} and the conclusions are still true for all $p>1$.
\end{rem}

%%%%%%%%%%%%%%%%%%%%%%%%%%%%%%%%%%%%%%%%%%%%%%%%%%%%%
\section{Diagonalization procedure}\label{diagprocedure}

We perform the Fourier transform of~\eqref{wavegeneral} with respect to $x$ obtaining
\begin{equation}
\label{eq:CPF}
\begin{cases}
\widehat{u}_{tt} +\langle \xi \rangle_{a, m}(t)^2\widehat{u}=0 \,,\\
\bigl(\widehat{u}(0,\xi),\widehat{u}_t(0,\xi)\bigr)=\bigl(\widehat{u_0}(\xi),\widehat{u_1}(\xi)\bigr)\,,
\end{cases}
\end{equation}
where $ \langle \xi \rangle_{a, m}(t) \doteq \left( |\xi|^2a(t)^2 + m(t)^2 \right)^{1/2}$.
 We put
\[ U=(i\langle \xi \rangle_{a, m}(t) \widehat{u},\widehat{u}_t)^T\,, \]
so  from \eqref{eq:CPF} we derive the system
\begin{equation}\label{eq:Uhyp}
\partial_tU=
\begin{pmatrix}
0 & 1 \\
1 & 0
\end{pmatrix}i\langle \xi \rangle_{a, m}(t)U + \frac{\partial_t \langle \xi \rangle_{a, m}(t)}{\langle \xi \rangle_{a, m}(t)} \begin{pmatrix}
1 & 0 \\
0 & 0
\end{pmatrix}U.
\end{equation}
Let $P$ be the (constant, unitary) diagonalizer of the principal part of~\eqref{eq:Uhyp}, given by
\[ P=\frac1{\sqrt{2}}\begin{pmatrix}
1 & 1 \\
-1 & 1
\end{pmatrix}\,, \qquad P^{-1}=\frac1{\sqrt{2}}\begin{pmatrix}
1 & -1 \\
1 & 1
\end{pmatrix}\,, \]
that is, if we put~$V(t,\xi)=P^{-1}U(t,\xi)$, then we get
\begin{equation}\label{eq:Vhyp}
\partial_tV=
\begin{pmatrix}
-1 & 0 \\
0 & 1
\end{pmatrix}i\langle \xi \rangle_{a, m}(t)V + R_1(t,\xi)\,V,
\end{equation}
where
\[ R_1(t,\xi) = \frac12 \frac{\partial_t \langle \xi \rangle_{a, m}(t)}{\langle \xi \rangle_{a, m}(t)} \begin{pmatrix} 1 & 1 \\
1 & 1 \end{pmatrix}. \]
We define the \emph{refined diagonalizer} which depends on the not diagonal entries of~$R_1(t,\xi)$:
\begin{equation}
\label{eq:Klambda}
K(t,\xi) \doteq I + K_1(t,\xi), \qquad K_1(t,\xi) \doteq
\frac1{4i} \frac{\partial_t \langle \xi \rangle_{a, m}(t)}{\langle \xi \rangle_{a, m}(t)^2} \begin{pmatrix}
0 & -1 \\
1 & 0
\end{pmatrix}\,.
\end{equation}
By using  Hypotheses~\ref{Hyp1} and~\ref{Hyp2}, if $\langle \xi \rangle_{a, m}(t) \eta(t)^{-1} \geq N$ we have
\begin{align}\label{eq:normK}
\left|\frac{\partial_t \langle \xi \rangle_{a, m}(t)}{\langle \xi \rangle_{a, m}(t)^2}\right| &=
\left|\frac{ a(t)a'(t)|\xi|^2+ m(t)m'(t)}{\langle \xi \rangle_{a, m}(t)^3}\right|\lesssim \frac{(a(t)^2|\xi|^2+ m(t)^2)\eta(t)}{\langle \xi \rangle_{a, m}(t)^3} \nonumber \\
 & \lesssim \frac{\eta(t)}{\langle \xi \rangle_{a, m}(t)} \leq \frac{C}N\,,
\end{align}
hence $\abs{\det K}\geq 1-C^2/N^2$. Therefore, $K(t,\xi)$ is uniformly regular and bounded for a sufficiently large $N$. We replace $V(t,\xi)=K(t,\xi)W(t,\xi)$ and we get
\begin{equation}\label{eq:W}
\partial_t W =
    \begin{pmatrix}
-1 & 0 \\
0 & 1
\end{pmatrix} i\langle \xi \rangle_{a, m}(t) W
    + \frac12 \frac{\partial_t \langle \xi \rangle_{a, m}(t)}{\langle \xi \rangle_{a, m}(t)}I W + R_2(t,\xi)W\,,
\end{equation}
 where the matrix $R_2$ is given by (see Lemma 5 in \cite{DAE2016})
\[R_2(t,\xi)= (\partial_t K_1+ K_1R_1) K^{-1}.\]
 Thanks again to Hypotheses~\ref{Hyp1} and~\ref{Hyp2},  the matrices~$R_2(t,\xi)$  satisfies the following estimate
\begin{equation}
\label{eq:Jestimate}
\|R_2(t,\xi)\|\lesssim  \frac{\eta(t)^2}{\langle \xi \rangle_{a, m}(t)}\,.
\end{equation}
Now let
\begin{equation}
\label{eq:D}
D(t,\xi) \doteq \diag \left(\exp \Bigl(-\int_{s}^t i\langle \xi \rangle_{a, m}(\tau)\,d\tau \Bigr) , \quad \exp \Bigl(\int_{s}^t i\langle \xi \rangle_{a, m}(\tau)\,d\tau \Bigr) \right)\,.
\end{equation}
We put $W(t,\xi)=\sqrt{\frac{\langle \xi \rangle_{a, m}(t)}{\langle \xi \rangle_{a, m}(s)}}\,D(t,\xi)Z(t,\xi)$ and we obtain %the following Cauchy problem
%in $Z_\hyp(N)$,
%
\begin{equation}
\label{eq:CPlambdaZ}
\begin{cases}
\partial_t Z = R_3(t,\xi)\,Z, \\
Z(s,\xi)=K^{-1}(s,\xi)P^{-1}\,U(s,\xi),
\end{cases}
\end{equation}
where the matrix $R_3(t,\xi)=D^{-1}(t,\xi)R_2(t,\xi)D(t,\xi)$ satisfies again~\eqref{eq:Jestimate}.

%%%%%%%%%%%%%%%%%%%%%%%%%%%%%%%%%%%%%%%%%%%%%%%%

\section{effective potential}\label{secteffect}

\begin{proof}(Theorem \ref{effectiveth})

We claim that
\begin{equation}\label{eq:claim}
\mathcal{E}(t,\xi)\leq C \gamma(t) \mathcal{E}_0(\xi)\,,
\end{equation}
uniformly with respect to~$\xi\in\R^n$, where~$\mathcal{E}(t,\xi)$ and~$\mathcal{E}_0(\xi)$ are  given by
\begin{align}
\label{eq:Etxi}
\mathcal{E}(t,\xi)\doteq |\widehat{u}_t(t,\xi)|^2+ (a(t)^2|\xi|^2 + \gamma(t)m(t))|\widehat{u}(t,\xi)|^2\,,\\
\label{eq:E0xi}
\mathcal{E}_0(\xi)\doteq |\widehat{u_1}(\xi)|^2+(1+|\xi|^2)|\widehat{u_0}(\xi)|^2\,.
\end{align}
Indeed, by integrating this inequality with respect to~$\xi$ and by Plancherel's Theorem, estimate~\eqref{effectiveRate} will follow from~\eqref{eq:claim}. In order to prove~\eqref{eq:claim}, for some constant~$N>0$, we divide the extended phase space~$[0,\infty)\times\R^n$ into the pseudo-differential zone $Z_\pd(N)$ and into the hyperbolic zone $Z_\hyp(N)$, defined by
\begin{align*}
Z_\pd(N)
    & =\{(t,\xi) \in [0,\infty) \times \R  : \langle \xi \rangle_{a, m}(t) \eta(t)^{-1} \leq N\},\\
Z_\hyp(N)
    & =\{(t,\xi) \in [0,\infty) \times \R : \langle \xi \rangle_{a, m}(t) \eta(t)^{-1} \geq N\}.
\end{align*}
By using the definition of $m$ we derive that
\[\langle \xi \rangle_{a, m}(t) \eta(t)^{-1}=(A(t)^2|\xi|^2+ \mu(t)^2)^{1/2}.\]
Therefore, thanks to the effective condition \eqref{effectivecondition},  $Z_\pd(N)$ is a compact subset of the extended phase space.
Then there exists a constant $T>0$ such that
 \[\mathcal{E}(t,\xi) \lesssim \mathcal{E}_0(\xi),\]
 for all $0\leq t\leq T$ and $(t,\xi)\in Z_\pd(N)$.

In $Z_\hyp(N)$ we use the calculations of Section \ref{diagprocedure}. Thanks to
$\frac{\eta(t)}{\mu(t)}\in L^1[0, \infty)$ we have
\begin{eqnarray*}
 \int_s^t \|R_3 (\tau,\xi)\|\,d\tau \lesssim \int_{s}^\infty \frac{\eta(\tau)^2}{\langle \xi \rangle_{a, m}(\tau)}\,d\tau
  \lesssim \int_{s}^\infty \frac{\eta(\tau)^2}{m(\tau)} \,d\tau=  \int_{s}^\infty  \frac{\eta(\tau)}{\mu(\tau)} \,d\tau
   \leq {C},
   \end{eqnarray*}
hence~$|Z(t,\xi)|\leq C|Z(s,\xi)|$ and, by using Liouville's formula, $|Z(t,\xi)|\geq C'|Z(s,\xi)|$. Therefore we have proved that in $Z_\hyp(N)$ it holds
\begin{equation}\label{eq:estimateZhyp}
C_1\,\frac{\langle \xi \rangle_{a, m}(t)}{\langle \xi \rangle_{a, m}(s)}\,|U(s,\xi)|^2 \leq |U(t,\xi)|^2 \leq C_2\,\frac{\langle \xi \rangle_{a, m}(t)}{\langle \xi \rangle_{a, m}(s)}\,|U(s,\xi)|^2\,.
\end{equation}
We remark that~\eqref{eq:estimateZhyp} is a two-sided estimate, that is, we have a precise description of the behavior of the energy in~$Z_\hyp(N)$.

Using again that $Z_\pd(N)$ is a compact subset, we have that $s\in [0,T]$ and
 $\langle \xi \rangle_{a, m}(s) \geq C \langle \xi \rangle_{a, m}(0)$.  Then, \eqref{eq:estimateZhyp} implies
 \[|\widehat{u}_t(t,\xi)|^2+ a(t)^2|\xi|^2 |\widehat{u}(t,\xi)|^2\lesssim \frac{\langle \xi \rangle_{a, m}(t)}{ (1+|\xi|^2)^{1/2}} \mathcal{E}_0(\xi)\lesssim
 (a(t) + m(t) ) \mathcal{E}_0(\xi),\]
 for all $t\geq 0$.
 Now, using again
 \eqref{eq:estimateZhyp}, for all $(t,\xi)\in Z_\hyp(N)$ we have that
\begin{align*}
 m(t)|\widehat{u}(t,\xi)|^2 &= m(t) \frac{\langle \xi \rangle_{a, m}^2(t)}{\langle \xi \rangle_{a, m}^2(t)}|\widehat{u}(t,\xi)|^2 \\
&\lesssim \frac{m(t)}{\langle \xi \rangle_{a, m}^2(t)}\frac{\langle \xi \rangle_{a, m}(t)}{\langle \xi \rangle_{a, m}(0)}
\left( \langle \xi \rangle_{a, m}^2(0)|\widehat{u}(0,\xi)|^2 + |\widehat{u_t}(0,\xi)|^2\right)\\
&\lesssim \frac{m(t)}{\langle \xi \rangle_{a, m}(t)}
\left( \langle \xi \rangle_{a, m}^2(0)|\widehat{u}(0,\xi)|^2 + |\widehat{u_t}(0,\xi)|^2\right) \\
&\lesssim
\left( \langle \xi \rangle_{a, m}^2(0)|\widehat{u}(0,\xi)|^2 + |\widehat{u_t}(0,\xi)|^2\right),
\end{align*}
thanks to
\[\frac{m(t)}{\langle \xi \rangle_{a, m}(t)}\lesssim \left(1+ \frac {|\xi|^2a(t)^2}{m(t)^2}\right)^{-1/2}\lesssim 1.\]
\end{proof}
%\begin{remark}
%On may relax the assumption $\frac{\eta(t)}{\mu(t)}\in L^1[0, \infty)$ in Theorem \ref{effectiveth}. Indeed,
%if we replace $Z_\hyp(N)$ by
%    \[Z_\hyp(N) =\{(t,\xi) \in [0,\infty) \times \R : |\xi|A(t) \geq N\},\]
%     then
%we may estimate
%\begin{eqnarray*}
% \int_s^t \|R_3 (\tau,\xi)\|\,d\tau \lesssim \int_{s}^\infty \frac{\eta(\tau)^2}{\langle \xi \rangle_{a, m}(\tau)}\,d\tau
%  \lesssim \frac1{|\xi|}\int_{s}^\infty \frac{a(\tau)}{A(\tau)^2} \,d\tau\lesssim \frac1{|\xi|A(s)}\leq N^{-1}
%      \end{eqnarray*}
%      and we derive the same estimates in $Z_\hyp(N)$. But now
%      \[Z_\pd(N)
%     =\{(t,\xi) \in [0,\infty) \times \R  : |\xi|A(t) \leq N\}\]
%     is no longer  a compact set. Here, to derive sharp results, one should add additional zones. (see for instance \cite{ERnew})
%\end{remark}
%%%%%%%%%%%%%%%%%%%%%%%%%%%%%%%%%%%%%%%%%%%%%%%%%

%%%%%%%%%%%%%%%%%%%%%%%%%%%%%%%%%%%%%%%%%%%%%%%%%%
\section{Non-effective potential}\label{noneffect}

 In order to get some feeling for the behavior of solutions to \eqref{wavegeneral} in the case of non-effective potential we can (see \cite{EKNR})  transform the  time-dependent potential to a time-dependent damping
and a new potential. If we introduce the change  of  variables given by $u(t,x)=\psi(t)v(t,x)$ the
 Cauchy problem \eqref{wavegeneral} takes the form
\begin{equation}\label{wavenoneffec2}
\begin{cases}
v_{tt} - a(t)^{2} \Delta v + 2\frac{\psi'}{\psi}(t) v_t + \left( \frac{\psi''}{\psi}(t) + m(t)^2\right)v=0,
\quad (t,x)\in (0,\infty)\times \R^n,
\\
(v(0,x),v_t(0,x))=(v_0(x),v_1(x)),
\quad x\in \R^n.
\end{cases}
\end{equation}
	Aiming to exclude contributions to the energy coming from  the time-dependent
	potential, thanks to the scattering assumption \eqref{ScatCond}, a sufficient condition is
	\begin{eqnarray}
	\frac{A(t)}{a(t)}\left( \frac{\psi''}{\psi}(t) + m(t)^2  \right) \in L^1.
	\end{eqnarray}
Under this assumption, we may use some ideas  developed in \cite{DAE2013}  to derive asymptotic properties of solutions to wave equations with
time-dependent non-effective dissipation (see also \cite{Wirth1}).

%%%%%%%%%%%%%%%%%%%%%%%%%%%%%%%%%%%%%%%%%%%%%%%%%%%%%%
Here we divide again the extended phase space $[0,\infty)\times\mathbb{R}^n$ into the {\it pseudo-differential zone} $Z_{pd}(N)$ and into the {\it hyperbolic zone}
$Z_{hyp}(N)$ which are defined by
\begin{eqnarray*}
	Z_{pd}(N)
	& =\{(t,\xi) \in [0,\infty) \times \R^n  :A(t)|\xi| \leq N\},\\
	Z_{hyp}(N)
	& =\{(t,\xi) \in [0,\infty) \times \R^n : A(t)|\xi| \geq N\}.
\end{eqnarray*}
The {\it separating curve} is given by
\[ \theta : (0,N]\to [0,\infty), \,\,\,\, \theta_{|\xi|}=A^{-1}(N/|\xi|).\]
We put also~$\theta_0=\infty$, and $\theta_{|\xi|}=0$ for any~$|\xi|\geq N$. The pair~$(t,\xi)$ from the extended phase space belongs to $Z_{pd}(N)$ (resp.
to~$Z_{hyp}(N)$) if and only if $t\leq\theta_{|\xi|}$ (resp.~$t\geq\theta_{|\xi|}$).

In the  $Z_{hyp}(N)$ we use the same energy and diagonalization procedure done in Section \ref{diagprocedure} to conclude
\begin{equation}\label{eq:estimateZhypnoneff}
C_1\,\frac{\langle \xi \rangle_{a, m}(t)}{\langle \xi \rangle_{a, m}(s)}\,|U(s,\xi)|^2 \leq |U(t,\xi)|^2 \leq C_2\,\frac{\langle \xi \rangle_{a, m}(t)}{\langle \xi \rangle_{a, m}(s)}\,|U(s,\xi)|^2\,,
\end{equation}
thanks to
\begin{eqnarray*}
 \int_s^t \|R_3 (\tau,\xi)\|\,d\tau \lesssim \int_{s}^\infty \frac{\eta(\tau)^2}{\langle \xi \rangle_{a, m}(\tau)}\,d\tau
  \lesssim \int_{s}^\infty \frac{a(\tau)}{|\xi|A(\tau)^2} \,d\tau\lesssim \frac1{|\xi|A(s)}
   \leq {C}.
   \end{eqnarray*}
But differently from the effective potential,  now $Z_{pd}(N)$ is no longer a compact subset,
so we have to apply a new strategy  to derive estimates in this zone:

\subsection*{Consideration in the pseudo-differential zone}
We will consider the following micro-energy in the  pseudo-differential zone
\[ V=\Big(\psi(t)\eta(t)\widehat{u}, \psi(t)\widehat{u}_t-\psi'(t)\widehat{u}\Big)^T, \,\,\, V_0(\xi)=\Big(\widehat{u}_0(\xi),
\widehat{u}_1(\xi)- \mu(0)^2\widehat{u}_0(\xi)\Big)^T. \] So we have
\begin{equation}\label{newsystem1}
\partial_t{V}(t,\xi)=\mathcal{A}(t,\xi){V}:=\left(\begin{array}{cc}
\frac{\eta'(t)}{\eta(t)} + 2 \frac{\psi'(t)}{\psi(t)}& \eta(t) \\
-\frac{1}{\eta(t)}\left(\frac{\psi''}{\psi}+\langle \xi \rangle_{a,m}(t)^2 \right)& 0
\end{array}\right){V}.
\end{equation}
We want to prove that the fundamental solution $E=E(t,s,\xi)$ to \eqref{newsystem1}, that is, the solution to
\[ \partial_tE=\mathcal{A}(t,\xi)E\,,\quad E(s,s,\xi)=I, \]
is bounded for all $t \in [0,\theta_{|\xi|}]$. If we put $E=(E_{ij})_{i,j=1,2}$, then we can write for
$j=1,2$ the following system of coupled integral equations of Volterra type:
\begin{eqnarray}
\label{eq:E1j}  E_{1j}(t,0,\xi) =\eta(t)\psi(t)^2 \left(\delta_{1j}+\int_0^t\frac1{\psi(\tau)^2}E_{2j}(\tau,0,\xi)d\tau\right),
\end{eqnarray}
\begin{align}
\label{eq:E2j}  E_{2j}(t,0,\xi) = \delta_{2j} - \int_0^t
\frac1{\eta(\tau)} \Big(\frac{\psi''}{\psi}(\tau)+\langle \xi \rangle_{a,m}(\tau)^2\Big)E_{1j}(\tau,0,\xi)d\tau.\hspace*{0.5cm}
\end{align}
By replacing \eqref{eq:E2j} into~\eqref{eq:E1j} and after
integration by parts we get
\begin{align}\label{E1j}
&&E_{1j}(t,0,\xi)=\eta(t)\psi(t)^2\Big(\delta_{1j}+\delta_{2j}\int_0^t\psi(\tau)^{-2}d\tau\Big) \nonumber\\
&&- \eta(t)\psi(t)^2 \int_0^t\frac1{\eta(\tau)}\Big(\frac{\psi''}{\psi}(\tau)+ \langle \xi \rangle_{a,m}(\tau)^2\Big)E_{1j}(\tau,0,\xi)\int_{\tau}^t\psi(s)^{-2}dsd\tau. %\hspace*{0.5cm}
\end{align}
 By using \eqref{neweq33} we conclude from \eqref{E1j} that
\[ |E_{1j}(t,0,\xi)|\leq C + C\int_0^t \frac1{\eta(\tau)}\left|\frac{\psi''(\tau)}{\psi(\tau)}
+\langle \xi \rangle_{a,m}(\tau)^2\right||E_{1j}(t,0,\xi)|d\tau.
\]
Applying Gronwall's type inequality  we conclude
\[ |E_{1j}(t,0,\xi)|\leq C \exp\Big(C\int_0^t \frac1{\eta(\tau)}\Big(\Big|\frac{\psi''(\tau)}{\psi(\tau)}+m(\tau)^2\Big|+|\xi|^2a(\tau)^{2}\Big)d\tau\Big).\]
In $Z_{pd}(N)$ we have  $A(t)|\xi|\leq N$. So, from the last estimate we get
\[ |E_{1j}(t,0,\xi)|\leq C \exp\Big(C\int_0^t \frac1{\eta(\tau)}\Big(\Big|\frac{\psi''(\tau)}{\psi(\tau)}+m(\tau)^2\Big|\Big)d\tau\Big).\]
Finally, by using again \eqref{neweq33} we get $\|E_{1j}(t,0,\xi)\|\leq C$. From the boundedness of $\|E_{1j}(t,0,\xi)\|$ we can estimate
$\|E_{2j}(t,0,\xi)\|\leq C$. Therefore, we proved
\begin{equation}\label{newconclusionpd}
\|V(t,\xi)\|\leq C\|V_0(\xi)\|\,\,\,\mbox{for all}\,\,\,t \in (0,\theta_{|\xi|}].
\end{equation}
%%%%%%%%%%%%%%%%%%%%%%%%%%%%%%%%%%%%%%%%%%%%%%%%%%%%%%%%%%%%%%%%%%%%%%%%%%%%%%
\begin{proof}(Theorem \ref{NoneffTeo})

 We claim that
\begin{equation}\label{apriorib}
|\widehat{u}_t(t,\xi)|^2 + a(t)^{2\ell}|\xi|^2|\widehat{u}(t,\xi)|^2\lesssim a(t)\left((1+|\xi|^2)|\widehat{u_0}(\xi)|^2 + |\widehat{u_1}(\xi)|^2\right)
\end{equation}
and
\begin{equation}\label{apriori}
|\widehat{u}(t,\xi)|^2\lesssim
\frac{1}{\eta(t)^2\psi(t)^2}\left(|\widehat{u_0}(\xi)|^2 +
\frac{|\widehat{u_1}(\xi)|^2}{1+|\xi|^2}\right),
\end{equation}
uniformly with respect to $\xi\in \R^n$. By integrating these inequalities with respect to~$\xi$ and by Plancherel's Theorem we have our desired estimate
\eqref{energyinequality}.

Let us first prove \eqref{apriorib}. By using Cauchy-Schwarz
inequality,
%($2ab\leq a^2/2 + 2b^2$)
\eqref{neweq32} and the considerations in the
pseudo-differential zone
%\textcolor{blue}{together with the increasing behavior of $\psi$}
we conclude for all $t\leq \theta_{|\xi|}$ the estimates
\begin{eqnarray*}
\frac{|V(t,\xi)|^2}{\psi(t)^2}&\geq& \eta(t)^2|\widehat{u}(t,\xi)|^2+ |\widehat{u}_t(t,\xi)|^2+\left|\frac{\psi'(t)}{\psi(t)}\right|^2|\widehat{u}(t,\xi)|^2-|\widehat{u}_t(t,\xi)
|\left|\frac{2\psi'(t)}{\psi(t)}\widehat{u}(t,\xi)\right|\\
&\geq& \eta(t)^2|\widehat{u}(t,\xi)|^2+ \frac{1}{2} |\widehat{u}_t(t,\xi)|^2-\left|\frac{\psi'(t)}{\psi(t)}\right|^2|\widehat{u}(t,\xi)|^2\\
&\geq& (1-c^2)\eta(t)^2|\widehat{u}(t,\xi)|^2+ \frac{1}{2} |\widehat{u}_t(t,\xi)|^2\\
&\geq&\frac{1-c^2}{N^2}a(t)^{2}|\xi|^2|\widehat{u}(t,\xi)|^2+ \frac{1}{2} |\widehat{u}_t(t,\xi)|^2.
\end{eqnarray*}
Therefore, by using  \eqref{newconclusionpd} we have for all $t\leq \theta_{|\xi|}$
\begin{equation}\label{additional}
|\widehat{u}_t(t,\xi)|^2+a(t)^{2}|\xi|^2|\widehat{u}(t,\xi)|^2\lesssim  |V(t,\xi)|^2\lesssim \frac{1}{\psi(t)^2}|{V_0}(\xi)|^2.
\end{equation}
In the  $Z_{hyp}(N)\cap \{|\xi|\geq N\}$, thanks to \eqref{eq:estimateZhypnoneff} we conclude
\[|\widehat{u}_t(t,\xi)|^2+a(t)^{2}|\xi|^2|\widehat{u}(t,\xi)|^2\lesssim  |U(t,\xi)|^2\lesssim a(t)\left((1+|\xi|^2)|\widehat{u_0}(\xi)|^2 + |\widehat{u_1}(\xi)|^2\right). \]

In the  $Z_{hyp}(N)\cap \{|\xi|\leq N\}$,
 we have to  glue the estimate
\eqref{newconclusionpd} with \eqref{eq:estimateZhypnoneff}. Putting $s=\theta_{|\xi|} $ and using  \eqref{eq:estimateZhypnoneff} we have
\[|\widehat{u}_t(t,\xi)|^2+a(t)^{2}|\xi|^2|\widehat{u}(t,\xi)|^2\lesssim
\frac{\langle \xi \rangle_{a, m}(t)}{\langle \xi \rangle_{a, m}(s)}\left(\langle \xi \rangle_{a, m}(s)^2|\widehat{u}(s,\xi)|^2+ |\widehat{u}_t(s,\xi)|^2  \right).  \]
Due to $A(s)|\xi|=N$  and $m(s)\leq |\xi|a(s)$ we have that
\[\frac{\langle \xi \rangle_{a, m}(t)}{\langle \xi \rangle_{a, m}(s)} \lesssim \eta(s)^{-1}\langle \xi \rangle_{a, m}(t)\lesssim a(t).\]
 Therefore, using that $\psi$ is non-decreasing and \eqref{additional} with $t=s$, we conclude
\[|\widehat{u}_t(t,\xi)|^2+a(t)^{2}|\xi|^2|\widehat{u}(t,\xi)|^2\lesssim a(t)\left(|\widehat{u_0}(\xi)|^2 + |\widehat{u_1}(\xi)|^2\right). \]

Now let us prove \eqref{apriori}. For $t\leq \theta_{|\xi|}$ we have
from \eqref{newconclusionpd} the estimate
$$|\widehat{u}(t,\xi)|^2\lesssim \frac{1}{\eta(t)^2\psi(t)^2}|{V_0}(\xi)|^2.$$
In order to estimate $|\widehat{u}(t,\xi)|^2$ in the hyperbolic zone
we split our considerations for $|\xi|\leq N$ and $|\xi|\geq N$. By
definition, $\theta_{|\xi|}=0$ for all $|\xi|\geq N$, and from
\eqref{eq:estimateZhypnoneff} we have
$$|\widehat{u}(t,\xi)|^2
\lesssim \frac1{a(t)} \left(|\widehat{u_0}(\xi)|^2 + \frac{|\widehat{u_1}(\xi)|^2}{|\xi|^2}\right) \,\,
\mbox{ for all}  \,\, |\xi|\geq N.$$  On the other hand, for $|\xi|\leq N$,
from \eqref{eq:estimateZhypnoneff} and \eqref{newconclusionpd} we conclude
\begin{eqnarray*}
&& |\widehat{u}(t,\xi)|^2\lesssim \frac1{\langle \xi \rangle_{a, m}(t)\langle \xi \rangle_{a, m}(s)}
\left(\langle \xi \rangle_{a, m}(s)^2|\widehat{u}(s,\xi)|^2+ |\widehat{u}_t(s,\xi)|^2  \right)\\
&& \lesssim  \frac{1}{\eta(t)\eta(s)\psi(s)^2}\left( |\widehat{u_0}(\xi)|^2 + |\widehat{u_1}(\xi)|^2   \right)\lesssim
\frac{1}{\eta(t)^2\psi(t)^2}\left( |\widehat{u_0}(\xi)|^2 + |\widehat{u_1}(\xi)|^2   \right)
\end{eqnarray*}
thanks to $\frac{1}{\eta(t)\psi(t)^2}$ be
increasing for large $t$. Using that $a\not\in L^1$, the  proof is completed.

\end{proof}

%%%%%%%%%%%%%%%%%%%%%%%%%%%%%%%%%%%%%%%%%%%%%%%%%%%%%
\section{Scattering Theory}\label{scattering}

In order to define our scattering operator, first we have to prove some a priori estimates for the fundamental solution
of a system associate to the Cauchy problem \eqref{wavegeneral}. To active this, we shall divide the extended phase
	space into two zones:
	\begin{align*}
	Z_{pd}(N) =& \left\{ (t,\xi) \in [0,\infty) \times \R^n : A(t)|\xi| \leq N \right\}, \\
	Z_{hyp}(N) =& \left\{(t,\xi) \in [0,\infty) \times \R^n : A(t)|\xi| \geq N \right\}.  \\
	\end{align*}
	The {\it separating curve} is given by
	\[ \theta : (0,N]\to [0,\infty), \,\,\,\, \theta_{|\xi|}=A^{-1}(N/|\xi|).\]
	We put also~$\theta_0=\infty$, and $\theta_{|\xi|}=0$ for any~$|\xi|\geq N$. The pair~$(t,\xi)$ from the extended phase space belongs to $Z_{pd}(N)$ (resp.
	to~$Z_{hyp}(N)$) if and only if $t\leq\theta_{|\xi|}$ (resp.~$t\geq\theta_{|\xi|}$).

	\textbf{Consideration in pseudo-differential zone.}

	We will consider the following micro-energy in the pseudo-diferential zone:\footnote{{In the definition of the micro-energy we will use $D_t\widehat{u}$, where
			$D_t=\frac{1}{i}\partial_t$.}}
	\begin{align} \label{MicroPseudo}
	U(t,\xi) = \frac{1}{\sqrt{a(t)}} (h(t,\xi)\widehat{u}, D_t \widehat{u})^T,
	\end{align}
	where $h$ is given be \eqref{hfunction}.
	Then
	\begin{align} \label{SisPseudo}
	D_t U = \widetilde{A}(t,\xi)U,
	\end{align}
	where
	\begin{align*}
	\widetilde{A}(t,\xi) = \begin{pmatrix}
	\frac{D_t h}{h}- \frac{D_t {a}}{2a} & h \\
	\frac{a(t)^2|\xi|^2 + m(t)^2}{h} & - \frac{D_t {a}}{2a}
	\end{pmatrix}.
	\end{align*}
	We want to prove that the fundamental solution $E=E(t,s,\xi)$ to \eqref{SisPseudo}, that is, the solution to
	$$ D_tE = \widetilde{A}(t,\xi) E \hspace{0.4cm} E(s,s,\xi)=I,   $$
	satisfies the estimate $\| E(t,s,\xi) \| \lesssim 1$. Indeed, in the pseudo-differential zone we shall only prove that the fundamental solution is bounded. If we put $E=(E_{ij})_{i,j=1,2}$, then we can write for $j=1,2$
	the following system of coupled integral equations of
	Volterra type:
	\begin{eqnarray}
	& E_{1j}(t,s,\xi) = \frac {h(t, \xi)}{\sqrt{a(t)}} \left( \frac{\delta_{1j}\sqrt{a(s)}}{h(s, \xi)} + i\int_s^t \sqrt{a(\tau)} E_{2j}(\tau,s,\xi)d\tau  \right) \label{EqPseudo1},  \\
	& E_{2j}(t,s,\xi) = \frac {1}{\sqrt{a(t)}} \left( \delta_{2j}\sqrt{a(s)}+ i \int_s^t \frac{a(\tau)^2|\xi|^2 + m(\tau)^2}{h(\tau, \xi)}\sqrt{a(\tau)}E_{1j}(\tau,s,\xi) d\tau \right). \label{EqPseudo2}
	\end{eqnarray}
	We may assume that $a(0)=1$.
	Replacing \eqref{EqPseudo2} in \eqref{EqPseudo1}, with $s=0$,
%integrating by parts we get
%	\begin{eqnarray}
%	E_{1j}(t,0,\xi) &=& \frac {h(t, \xi)}{\sqrt{a(t)}}\left( \delta_{1j}   +i   \int_s^t\delta_{2j} d\tau\right) \nonumber  \\
%	&-& \frac {h(t, \xi)}{\sqrt{a(t)}}\int_0^t \int_0^{\tau} \frac{a(\sigma)^2|\xi|^2 + m(\sigma)^2}{h(\sigma, \xi)}\sqrt{a(\sigma)}E_{1j}(\tau,0,\xi)d\sigma d\tau . \label{EqScat}
%	\end{eqnarray}
	and using that $a(t)$ is an increasing function  we have
	\begin{eqnarray*}
		|E_{1j}(t,0,\xi)| &\leq & \frac {h(t, \xi)}{\sqrt{a(t)}}\left( \delta_{1j}   +  t \delta_{2j}\right) \nonumber  \\
		&+&  \frac {h(t, \xi)}{\sqrt{a(t)}} \int_0^t \sqrt{a(\tau)}\int_0^{\tau} \frac{a(\sigma)^2|\xi|^2 + m(\sigma)^2}{h(\sigma, \xi)}|E_{1j}(\sigma,0,\xi)|d\sigma d\tau .
	\end{eqnarray*}
	Integration by part  yields
	\begin{eqnarray*}
		|E_{1j}(t,0,\xi)| &\leq & \frac {h(t, \xi)}{\sqrt{a(t)}}\left( \delta_{1j}   +   t\delta_{2j} \right) \nonumber  \\
		&+&  \frac {h(t, \xi)}{\sqrt{a(t)}} \int_0^t \frac{a(\tau)^2|\xi|^2 + m(\tau)^2}{h(\tau, \xi)}|E_{1j}(\tau,0,\xi)| \left(\int_{\tau}^t \sqrt{a(\sigma)}d\sigma\right) d\tau .
	\end{eqnarray*}
	From the definition of  the pseudo-differential zone we can get
	\[
	\frac{A(t)}{a(t)}{h}(t,\xi) \sim  1,
	\]
	and thanks to hypothesis  \eqref{additionalscat} we have that
	\[\frac {h(t, \xi)}{\sqrt{a(t)}}\left( 1+t+ \int_{0}^t \sqrt{a(\sigma)}d\sigma\right)\leq C.\]
	Hence, applying Gronwal's inequality, condition \eqref{ScatCond} and from the definition of the  pseudo-differential zone we can get
	\[
	|E_{1j}(t,0,\xi)| \leq  \exp\left( \int_0^t \frac{a(\sigma)^2|\xi|^2 + m(\sigma)^2}{h(\sigma, \xi)} d\tau\right)\leq C.
	\]
	Now, using again  \eqref{ScatCond}  and that $|E_{1j}(t,0,\xi)|$  is bounded in \eqref{EqPseudo2}, we also conclude that $|E_{2j}(t,0,\xi)| \leq C$.
	\
	
	\textbf{Consideration in hyperbolic zone.}
	
	Define the micro-energy
	$$ U_W(t,\xi)= (a(t)|\xi|\widehat{u}, D_t\widehat{u})^T. $$
	Then $U_W $ satisfies
	\begin{align*}
	D_tU_W = A(t,\xi) U_W,
	\end{align*}
	with
	$$ A(t,\xi) = \begin{pmatrix}
	\frac{D_t a(t)}{a(t)} & a(t)|\xi| \\
	a(t)|\xi| + \frac{m(t)^2}{a(t)|\xi|} & 0
	\end{pmatrix} $$
	Let
	$$ M = \begin{pmatrix}
	1 & -1 \\
	1 & 1
	\end{pmatrix} \,\,\, \mbox{and} \,\,\, M^{-1} = \frac{1}{2}\begin{pmatrix}
	1 & 1 \\
	-1 & 1
	\end{pmatrix}.$$
	If we define $U^{(0)} = M^{-1}U_W,$ then
	\begin{align*}
	D_t U^{(0)}  = \left( D(t,\xi) + R_a(t) + R_{a,m}(t)   \right)U^{(0)},
	\end{align*}
	with
	\begin{align*}
	D(t,\xi) = \begin{pmatrix}
	a(t)|\xi|& 0 \\
	0 & -a(t)|\xi|
	\end{pmatrix}, \,\,\, R_a(t) = \frac{1}{2} \frac{D_ta(t)}{a(t)}\begin{pmatrix}
	1 & -1 \\
	-1 & 1
	\end{pmatrix}
	\end{align*}
	and
	\begin{align*}
	R_{a,m}(t, \xi) =  \frac{1}{2}  \frac{m(t)^2}{a(t)|\xi|}\begin{pmatrix}
	1 & -1 \\
	-1 & 1
	\end{pmatrix}.
	\end{align*}
	Let $E_a = E_a(t,s,\xi)$ be the fundamental solution of the operator $D_t - D(t,\xi) - R_a(t)$, that is, $E_a$ satisfies the Cauchy problem
	$$ \left(D_t - D(t,\xi) - R_a(t)\right) E_a(t,s,\xi) =0, \,\,\, E(s,s,\xi) = I.  $$
	It is well known (see \cite{ReisBui}) that
	\[ \| E_a(t,s,\xi) \|_{L^\infty(\R^n_\xi)} \lesssim \frac{\sqrt{a(t)}}{\sqrt{a(s)}}. \]
	Moreover, using Liouville's formula, we arrive at
	\[ \det E_a(t,s,\xi) = \exp\left( i \int_s^t tr(D(\tau, \xi)+R_a(\tau))d\tau \right) = \frac{a(t)}{a(s)}.\]
	Hence,
	\[ \| E_a^{-1}(t,s,\xi) \|_{L^\infty(\R^n_\xi)} \lesssim \frac{\sqrt{a(s)}}{\sqrt{a(t)}}. \]
	Now the goal is to construct the fundamental solution to the operator $D_t-D(t,\xi) -R_a(t) -R_{a,m}(t).$ For this purpose let us introduce
	$$  P(t,s,\xi) = E_a(t,s,\xi)^{-1} R_{a,m}(t, \xi)E_a(t,s,\xi).  $$
	Applying Peano-Backer formula we have that
	\begin{align}\label{SolRb}
	Q_{a,m}(t,s,\xi)= I + \sum_{k=1}^\infty i^k\int_s^tP(t_1,s,\xi)
	\int_s^{t_1}P(t_2,s,\xi) \cdots \int_s^{t_{k-1}}P(t_k,s,\xi)dt_k \cdots dt_2 dt_1
	\end{align}
	is the solution to the Cauchy problem
	$$ D_t Q_{a,m}(t,s,\xi) = P(t,s,\xi)Q_{a,m}(t,s,\xi), \,\,\, Q_{a,m}(s,s,\xi) = I.  $$
	Therefore, if $E_{a,m}(t,s,\xi)=E_a(t,s,\xi)Q_{a,m}(t,s,\xi)$ follows that
	\begin{align}
	D_t E_{a,m}(t,s,\xi) = (D(t,\xi)+R_a(t)+ R_{a,m}(t, \xi)) E_{a,m}(t,s,\xi), \hspace{0.3cm} E_{a,m}(s,s,\xi)=I.
	\end{align}
	From the definition of $P(t,s,\xi)$ we can derive that
	$$ \| P(t,s,\xi) \|_{L^\infty} \leq \|  R_{a,m}(t, \xi) \|. $$
	So, using the definition of the hyperbolic zone and the hypothesis \eqref{ScatCond} we arrive at
	\[  \| Q_{a,m}(t,s,\xi) \|_{L^\infty} \leq \exp \left( \int_s^t \|P(\tau, s, \xi)  \|_{L^\infty} d\tau \right)
\leq \exp \left( \int_s^t \frac{A(\tau)m(\tau)^2}{a(\tau)} d\tau \right) \leq C. \]
	Consequently,
	\[ \| E_{a,m}(t,s,\xi) \|_{L^\infty} \lesssim \|E_a(t,s,\xi) \|_{L^\infty}\|Q_{a,m}(t,s,\xi) \|_{L^\infty} \lesssim \frac{\sqrt{a(t)}}{\sqrt{a(s)}}. \]
	Now, let us introduce
	\[H(t,\xi) :=\begin{pmatrix} \frac{h(t,\xi)}{|\xi|a(t)} & 0 \\
	0 & 1 \end{pmatrix}.\]
	It is clear that in the hyperbolic zone we have $\frac{h(t,\xi)}{|\xi|a(t)}\approx C$. Then the inverse matrix $H^{-1}$ exists and $\|
	H(t,\xi)\|, \| H^{-1}(t,\xi)\| \approx C$ for all $t\geq \theta_{|\xi|}$.
	
	Thanks to
	\[U(t,\xi)=\frac{1}{\sqrt{a(t)}}HU_W=\frac{1}{\sqrt{a(t)}}( h(t,\xi)\widehat{u},D_t\widehat{u})^T,\]
	where  $U$ is  defined  in \eqref{MicroPseudo}, we conclude from the previous calculation in the hyperbolic zone that
	\[ |U(t,\xi)|=|\frac{\sqrt{a(s)}}{\sqrt{a(t)}}H(t,\xi)ME_a(t,s,\xi)Q_{a,m}(t,s,\xi)M^{-1}H^{-1}(s,\xi) U(s,\xi)| \leq C |U(s,\xi)|. \]
	
	%%%%%%%%%%%%%%%%%%%%%%%%%%%%%%%%%%%%%%
	
	\begin{flushleft}
		
	\end{flushleft}%%%%%%%%%%%%%%%%%%%%%%%%%%%%%%%%%%%%%%%%%%%%%%%%

\begin{proof}( Theorem \ref{scat})
	
	With $s=\theta_{|\xi|}$ and the notation introduced in the pseudo-differential zone and in the hyperbolic zone we define
	$$\mathcal{E}(t, s,\xi)=\left\{
	\begin{array}{ccc}
	& E(t,0,\xi), \ 0\leq t\leq \theta_{|\xi|}, \hspace{9cm}\\
	&\frac{\sqrt{a(s)}}{\sqrt{a(t)}} H(t,\xi)ME_a(t,s,\xi)Q_{a,m}(t,s,\xi)M^{-1}H^{-1}(s,\xi)E(s,0,\xi) , \ t\geq \theta_{|\xi|}. \hspace{3cm}
	\end{array}
	\right.$$
	We have proved that $\|\mathcal{E}(t, s,\xi)\|\leq C$ for all $t, \xi$.
	
	Note that the matrix functions  $\mathcal{E}(t,s,\xi)$ and  $\frac{\sqrt{a(s)}}{\sqrt{a(t)}}ME_a(t,s,\xi)M^{-1}$
	generates a Fourier multiplier to the operators
\[ S(t,s,D): \frac{1}{\sqrt{a(s)}}( \langle D(s)\rangle u(s),D_tu(s))^T  \mapsto \frac{1}{\sqrt{a(t)}}( \langle D(t)\rangle u(t),D_tu(t))^T,  \]
	\[ S_1(t,s,D): \frac{1}{\sqrt{a(s)}}(a(s)|D|v(s),D_tv(s))^T  \mapsto \frac{1}{\sqrt{a(t)}}(a(t)|D|v(t),D_tv(t))^T,  \]
		for the solutions $u$ and $v$ to the Cauchy problems \eqref{wavegeneral} and  \eqref{CauchyProbSpeeWave}, respectively.
	Therefore we shall prove that the limit
	$$ W_+(D) = \lim_{t\rightarrow \infty} S_1^{-1}(t,0,D) S(t,0,D)$$
	exists in $E$. To study the operator $S_1^{-1}(t,0,D) S(t,0,D)$ it is sufficient to study in the phase space the  bounded multiplier
	\[ ME_a^{-1}(t,0,\xi)M^{-1} H(t,\xi)ME_a(t,s,\xi)Q_{a,m}(t,s,\xi)M^{-1}H^{-1}(s,\xi)E(s,0,\xi) .\]
	Let us prove that the limit $W_+(\xi)$ exists for all $|\xi|\geq c>0$.  Thanks to $E_a^{-1}(t,0,\xi)=E_a(0,t,\xi)$, $E_a(0,t,\xi)E_a(t,s,\xi)=E_a(0,s,\xi)$ and  $\lim_{t\rightarrow \infty} H(t,\xi)=I$, we get that
	\[\lim_{t\rightarrow \infty} E_a^{-1}(t,0,\xi)M^{-1} H(t,\xi)ME_a(t,s,\xi)=E_a(0,s,\xi).\]
	Thus we shall investigate the limit
	\begin{align}\label{LimOpe}
	\lim_{t \rightarrow \infty} Q_{a,m}(t,\theta_{|\xi|},\xi).
	\end{align}
	We obtain from the Peano-Backer formula that
	\begin{eqnarray*}
		Q_{a,m}(t,0,\xi) -Q_{a,m}(s,0,\xi) =\\
		\sum_{k=1}^\infty i^k \int_s^t P(t_1,0,\xi)
		\int_0^{t_1} P(t_2,0,\xi) \cdots \int_0^{t_{k-1}} P(t_k,0,\xi) dt_k \cdots dt_2 dt_1.
	\end{eqnarray*}
	Therefore,
	\begin{eqnarray*}
		&\|Q_{a,m}(t,\theta_{|\xi|},\xi) -Q_{a,m}(s,\theta_{|\xi|},\xi)\|_{L^\infty} \leq \sum_{k=1}^\infty  \int_s^t \|P(t_1,\theta_{|\xi|},\xi))\|_{L^\infty} \\
		&\times \frac{1}{(k-1)!} \left(\int_{\theta_{|\xi|}}^{t_1} \|P(\tau,\theta_{|\xi|},\xi))\|_{L^\infty}\right)^{k-1}dt_1\\
		&\leq \int_s^t \|P(t_1,\theta_{|\xi|},\xi))\|_{L^\infty} \sum_{k=0}^\infty \frac{1}{k!} \left(\int_{\theta_{|\xi|}}^{t_1}\|P(\tau,\theta_{|\xi|},\xi))\|_{L^\infty}\right)^{k}dt_1 \\
		&\leq \int_s^t \| R_{a,m}(t_1, \xi) \| \exp \left( \int_{\theta_{|\xi|}}^{t_1} \| R_{a,m}(\tau, \xi) \|\d\tau\right)dt_1.
	\end{eqnarray*}
	Our assumption \eqref{ScatCond} implies that in the hyperbolic zone $R_{a,m}(\cdot, \xi) \in L^1$. Thus,
	\[ \|Q_{a,m}(t,\theta_{|\xi|},\xi) -Q_{a,m}(s,\theta_{|\xi|},\xi)\|_{L^\infty} \]
becomes arbitrarily small for sufficiently large times $s,t$.
	 So the limit \eqref{LimOpe} exists uniformly in $\xi$ for $|\xi|\geq c>0$  and thanks to the Banach-Steinhaus's theorem we can define the  operator $W_+$.

The conclusion of the theorem follows thanks to
\begin{eqnarray*}
&\frac{1}{\sqrt{a(t)}} \Big((a(t)\nabla v(t,\cdot),v_t(t,\cdot))-(\langle D(t)
	\rangle u(t,\cdot),u_t(t,\cdot))\Big)\\
&=S_1(t,0,D)\Big (S_1^{-1}(t,0,D)S(t,0,D) - W_+(D) \Big)  (\langle D(0)\rangle u_0,u_1),
\end{eqnarray*}
and for the decay rate we only use that
	\begin{eqnarray*}
		&\|  Q_{a,m}(t,0,\xi) - Q_{a,m}(\infty,0,\xi) \|_{L^\infty}  \lesssim  \int_t^\infty \| R_{a,m}(t_1) \| \exp \left( \int_0^{t_1} \| R_{a,m}(\tau) \|\d\tau\right)dt_1
		\\
		&\lesssim \int_t^\infty \frac{A(\tau)}{a(\tau)}m(\tau)^2 d\tau.
	\end{eqnarray*}
	The proof is completed.
\end{proof}

%%%%%%%%%%%%%%%%%%%%%%%%%%%%%%%%%%%%%%%%%%%%%%%%

\section{Scale invariant models}\label{scaleinvariant}
\label{SecOptiCritMo}
Let us consider the Cauchy problem for a class of scale invariant models with time-dependent mass and speed of propagation
\begin{eqnarray}\label{CauchyProb1}
u_{tt} - a(t)^{2} \Delta u  + m(t)^2 u
=0,\,\,\,u(0,x)=u_0(x),\,\,\,u_t(0,x)=u_1(x),
\end{eqnarray}
with $m(t)=\mu\frac{a(t)}{A(t)}$ and the function $a\notin L^1$ is given by $a(t)= A(0)^{-\alpha}A(t)^{\alpha}$ for some constant $\alpha \in \mathbb{R}$, i.e.,
\[\frac{a'(t)}{a(t)} = \alpha\frac{a(t)}{A(t)}, \qquad A(t)  \doteq A(0)+ \int_0^t a(\tau)d\tau.\]
Assume that $A(0)\in (0,1]$. Applying the change of variable
\begin{align}
\label{ChangOfVariCM}
v(\tau, x)=u(t,x), \qquad \tau+1=A(t), \qquad \tau_0\doteq A(0)-1\in (-1, 0]
\end{align}
the Cauchy problem  \eqref{CauchyProb1} takes the form
\begin{align} \label{CauchPrAfCV}
v_{\tau\tau} -  \Delta v  + \frac{\alpha}{1+\tau} v_{\tau} + \frac{\mu^2}{(1+\tau)^2} v=0,\,\,\,v(\tau_0,x)=u_0(x),\,\,\,v_{\tau}(\tau_0,x)=u_1(x).
\end{align}
Now, the number $\delta\doteq(\alpha-1)^2-4\mu^2$ plays a fundamental role.

\begin{itemize}
\item If $\delta < 0$, applying the change of variable
$v(\tau, x)=(1+\tau)^{-\alpha/2}w(\tau,x)$ (see \cite{NPR})
we get the equation
\[w_{\tau\tau} -  \Delta w  + \frac{\sigma}{(1+\tau)^2} w=0,\,\,\,w(\tau_0,x)=w_0(x),\,\,\,w_{\tau}(\tau_0,x)=w_1(x),\]
where $\sigma=\frac{\alpha}{2}-\frac{\alpha^2}{4} + \mu^2=\frac{1-\delta}{4}\geq 1/4$ and
\[w_0(x)\doteq (1+\tau_0)^{\alpha/2}u_0(x), \qquad  w_1(x)\doteq  \frac{\alpha}{2} (1+\tau_0)^{-1}u_0(x)+
(1+\tau_0)^{\alpha/2}u_1(x).\]
 Applying a result from \cite{B1} we get
that the solution to \eqref{CauchyProb1} satisfies
\begin{eqnarray}\label{energyestimates}
E(u)(t) \lesssim a(t) E(u)(0), \qquad \forall t\geq 0,
\end{eqnarray}
where
\[E(u) (t)\doteq \frac{1}{2} \left(\|u_t(t,\cdot) \|_{L^2}^2 + a(t)^2\|\triangledown_x u(t,\cdot) \|_{L^2}^2 + m(t)a(t)\|u(t,\cdot) \|_{L^2}^2 \right).\]
In particular we get
\begin{align} \label{EstForThSolCC}
\|u(t,\cdot) \|_{L^2} \lesssim \frac1{\sqrt{m(t)}} E(u)(0) \sim \frac{\sqrt{A(t)}}{\sqrt{a(t)}} E(u)(0) , \qquad \forall t\geq 0.
\end{align}

%\begin{remark}
%If $\alpha>1$, then $m(t)=\mu A(t)^{\alpha -1}$ is an increasing function and $\|u(t,\cdot) \|_{L^2}$ decay.
%\end{remark}

\item If $\delta \geq0$,  applying the change of variable (see \cite{NPR})
\[v(\tau, x)=(1+\tau)^{\sigma}w(\tau,x), \qquad \sigma\doteq\frac{1-\alpha}{2} + \frac{\sqrt \delta}{2},\]
we get the equation
\[w_{\tau\tau} -  \Delta w  + \frac{1+\sqrt \delta}{(1+\tau)} w_{\tau}=0,\,\,\,w(\tau_0,x)=w_0(x),\,\,\,w_{\tau}(\tau_0,x)=w_1(x),\]
where
\[w_0(x)\doteq (1+\tau_0)^{-\sigma}u_0(x), \qquad  w_1(x)\doteq   (1+\tau_0)^{-\sigma}u_1(x)- \sigma (1+\tau_0)^{ -1}a(\tau_0)^{-1}u_0(x) .\]
If $(u_0,u_1) \in H^1 \times L^2$, applying  results from \cite{W} we get
\begin{eqnarray*}
\|w(\tau,\cdot) \|_{L^2}\lesssim \left\{ \begin{array}{cr}
1 , &\delta>0,\\
\ln (e+\tau),& \delta=0, \end{array} \right.\end{eqnarray*}
and
\begin{eqnarray*}
\|w_{\tau}(\tau,\cdot) \|_{L^2}+ \| \nabla w(\tau,\cdot) \|_{L^2}\lesssim \left\{ \begin{array}{cr}
(1+\tau)^{-\frac{(1+\sqrt \delta)}{2}} , &\delta \in [0,1),\\
(1+\tau)^{-1},& \delta\geq 1. \end{array} \right.\end{eqnarray*}

Therefore, using that the solution to \eqref{CauchyProb1} satisfies $u(t, x)=(1+\tau)^{\sigma}w(\tau,x)$ we conclude
\begin{eqnarray} \label{eq:w}
\|u(t,\cdot) \|_{L^2}\lesssim \left\{ \begin{array}{cr}
a(t)^{\frac{1-\alpha + \sqrt \delta}{2\alpha}} , &\delta>0,\\
a(t)^{\frac{1-\alpha }{2\alpha}}\ln a(t),& \delta=0, \end{array} \right.\end{eqnarray}
and
\begin{eqnarray*}
a(t)\| \nabla u(t,\cdot) \|_{L^2}\lesssim \left\{ \begin{array}{cr}
a(t)^{\frac{2(\sigma +\alpha)-1-\sqrt \delta}{2\alpha}} , &\delta \in [0,1),\\
a(t)^{\frac{\sigma -1+\alpha}{\alpha}},& \delta\geq 1.\end{array} \right.\end{eqnarray*}
%By using that $\alpha>0 (< 0)$ for increasing (decreasing) $a(t)$
We may derive that same estimate for $\|u_{t}(t,\cdot) \|_{L^2}$
and thanks to $\frac{2(\sigma +\alpha)-1-\sqrt \delta}{\alpha}=1$ we conclude
\begin{eqnarray} \label{eq:w1}
\|u_{t}(t,\cdot) \|_{L^2}+ a(t)\| \nabla u(t,\cdot) \|_{L^2}\lesssim \left\{ \begin{array}{cr}
\sqrt{a(t)}, &\delta \in [0,1),\\
a(t)^{\frac{\alpha-1+\sqrt\delta}{2\alpha}},& \delta\geq 1.\end{array} \right.\end{eqnarray}
\end{itemize}
\begin{exam}\label{exexp}(Exponential speed of propagation)\\
Consider the  Cauchy problem
\[
u_{tt} - e^{2t} \Delta u  + \mu^2 u
=0,\,\,\,u(0,x)=u_0(x),\,\,\,u_t(0,x)=u_1(x),
\]
 i.e., model \eqref{CauchyProb1} with $\mu>0$ and $\alpha=1$.  Then $\delta= -4\mu^2<0$ and thanks to \eqref{energyestimates} we conclude that
\[E(u) (t) := \frac{1}{2} \left(\|u_t(t,\cdot) \|_{L^2}^2 + e^{2t}\|\triangledown_x u(t,\cdot) \|_{L^2}^2 +\mu e^{t}\|u(t,\cdot) \|_{L^2}^2 \right)
\lesssim  e^{t}E(u)(0).\]
\end{exam}
\begin{exam}(Polynomial speed of propagation)\label{polinomial}\\
Consider the  Cauchy problem
\[
u_{tt} - (1+t)^{2\ell} \Delta u  + \frac{{\tilde \mu}^2}{(1+t)^2} u
=0,\,\,\,u(0,x)=u_0(x),\,\,\,u_t(0,x)=u_1(x),
\]
where $\ell>-1$ and $\tilde{\mu}>0$. This model can be written in the form \eqref{CauchyProb1} with $\alpha= \frac{\ell}{\ell+1}$ and $\tilde{\mu}=\mu(\ell+1)$.
In this case
\[ \delta=\frac{ 1-4 {\tilde \mu}^2}{(\ell+1)^2}.\]
	If $(u_0,u_1) \in H^1 \times L^2,$ thanks to \eqref{energyestimates}, \eqref{eq:w} and \eqref{eq:w1} we have the following estimates
\begin{eqnarray*} \label{eq: p(t)}
\|u(t,\cdot) \|_{L^2}^2 \lesssim \left\{ \begin{array}{cr}
(1+t), & \delta<0,\\
(1+t)(\ln (e+t))^2,&  \delta=0,\\
(1+t)^{1+\sqrt{1-4{\tilde \mu}^2}} ,&   \delta>0, \end{array} \right.\end{eqnarray*}
and
\begin{eqnarray*} \label{eq: p(t)}
\|u_t(t,\cdot) \|_{L^2}^2 +(1+t)^{2\ell}\|\nabla u(t,\cdot) \|_{L^2}^2 \lesssim \left\{ \begin{array}{cr}
(1+t)^{\ell}, & \delta<1,\\
(1+t)^{-1+\sqrt{1-4{\tilde \mu}^2}},& \delta \geq 1. \end{array} \right.\end{eqnarray*}
	\end{exam}
\begin{remark}
In Example \ref{polinomial},  the case $\delta\geq 1$ correspond to $-1<\ell<0$ satisfying
$\ell+1\leq \sqrt {1-4{\tilde \mu}^2}$. In particular,  this condition is equivalent to say that for negative $\ell$,
the last decay for the kinetic and elastic  energies  is worst than the first one in the case  $\delta<1$.
\end{remark}

%%%%%%%%%%%%%%%%%%%%%%%%%%%%%%%%%%%
%%%%%%%%%%%%%%%%%%%%%%%%%%%%%%%%%%%
%%%%%%%%%%%%%%%%%%%%%%%%%%%%%%%%

\subsection{$L^q-L^2$ estimates, $q \in [1,2)$ }\label{generalestimates}

\

In this section we show that additional regularity  $L^q$, with $q \in [1,2)$, may improve the estimates for the solution and its derivatives.
As discussed in the Section \ref{SecOptiCritMo}, applying the change of variable
\eqref{ChangOfVariCM} we arrive in the Cauchy problem \eqref{CauchPrAfCV}. Depending on the signal of $\delta= (\alpha-1)^2-4\mu^2$ we have two situations:

\begin{enumerate}
		\item If $\delta \leq 0,$ $(u_0,u_1) \in \mathcal{D}_q$, then applying  Theorem 4.3 from \cite{NPR}  the solution
		$u $
		for the Cauchy problem  \eqref{CauchyProb1}
		satisfies the following estimates:
		\[
		\|  u(t,\cdot) \|_{L^2} \lesssim \frac{1}{\sqrt{a(t)}}( 1+\ln (a(t)^{\frac{1}{\alpha}}  ))^\gamma d(t) \big( \| u_0\|_{H^1\cap L^q} + \| u_1 \|_{L^2\cap L^q}\big),
		\]
		for all $t\geq  0$, where $\gamma =1$ if $\delta =0$, $\gamma =0$ if $\delta < 0$ and
		%\begin{align*}
		% q_{0}(t,s)= \left\{ \begin{array}{cr}
		% 1+\ln \Big( \frac{1+t}{1+s} \Big)  &\mbox{\,\,\, for \,\,\,}  n>\frac{m}{2-m},\\
		% \Big(\ln\big(\frac{1+t}{1+s}\big) \Big)^{\frac{2-m}{2}} \Big( 1+\ln \Big( \frac{1+t}{1+s} \Big) \Big)   &\mbox{\,\,\, for \,\,\,} n=\frac{m}{2-m}, \end{array} \right.
		%\end{align*}
		%and
		\begin{align*}
		d(t)= \left\{ \begin{array}{cr}
		1  & \mbox{\,\,\, for \,\,\,} n>\frac{q}{2-q},\\
		( \ln(a(t)^{\frac{1}{\alpha}} )^{\frac{2-q}{2q}}  &  \mbox{\,\,\, for \,\,\,} n=\frac{q}{2-q}. \end{array} \right.
		\end{align*}
			 We notice that additional regularity in the initial data improves the estimates for the potential energy, see \eqref{EstForThSolCC}  for the case $q=2$.
			
		\item If $\delta>0,$ $(u_0,u_1) \in \mathcal{D}_q^{\kappa-1}$ and $\kappa\in[0,1]$ then  applying Theorem 4.7 from \cite{NPR} the solution
		$u $
		for the Cauchy problem  \eqref{CauchyProb1}
		satisfies the following estimates:
		\begin{equation*}\label{Lin Estim Eff Diss} ||u(t,\cdot)||_{\dot{H}^\kappa}\lesssim ||(u_0,u_1)||_{\mathcal{D}_q^{\kappa-1}} \begin{cases}a(t)^{\frac{1}{\alpha}\left(-\kappa- \frac{2-q}{2q}n +\frac{1-\alpha+\sqrt{\delta}}{2}\right)}, &  1+\sqrt{\delta} >\frac{2-q}{q}n+2\kappa, \\
		a(t)^{-\frac{1}{2}}\left(1+(\log(a(t)^{\frac1{\alpha}}))^{\frac{2-q}{q}}\right), &   1+\sqrt{\delta} =\frac{2-q}{q}n+2\kappa, \\
		a(t)^{-\frac{1}{2}}, &   1+\sqrt{\delta} <\frac{2-q}{q}n+2\kappa. \end{cases}
		\end{equation*}
Moreover, $a(t)^{-1}||\partial_t u(t,\cdot)||_{L^2(\mathbb{R}^n)}$  satisfies the same decay estimates as \\$||\nabla u(t,\cdot)||_{L^2(\mathbb{R}^n)}$ which are obtained from \eqref{Lin Estim Eff Diss} after taking $\kappa=1$.

In this case we notice that additional regularity on the initial data only improves the estimates for the potential energy if $a(t)$ is increasing, see  \eqref{eq:w} for the case $q=2$.

%\begin{rem}
%If $a(t)$ is increasing (decreasing), then the \textcolor{red}{ decay rate} for the potential energy is better for small (large) values of $\delta$.
%\end{rem}
\end{enumerate}

\subsection{Proof of Theorem \ref{nonlinear}}

\

 According to Duhamel's principle, a solution to \eqref{CauchyProbN2} satisfies the non-linear integral equation
\[  u(t,x)=K_0(t,0, x) \ast_{(x)} u_0(x) + K_1(t,0, x) \ast_{(x)} u_1(x) + \int_0^t K_1(t,s,x) \ast_{(x)}|u(s,x)|^{p}\,ds, \]
where $K_j(t,0, x) \ast_{(x)} u_j(x)$, $j=0,1$, are the solutions to the corresponding linear Cauchy
problem
\begin{equation}\label{sitterwave}
u_{tt} - e^{2t}\Delta u + m^2 u=0,
  \quad
u(0,x)=\delta_{0j}u_0(x), \quad u_t(0,x)=\delta_{1j}u_1(x),
  \end{equation}
with $\delta_{kj}=1$  for $k=j$, and zero otherwise. The term $K_1(t,s,x) \ast_{(x)}f(s,x)$ is the  solution  of the parameter-dependent Cauchy problem
\begin{equation}\label{linearassociate}
u_{tt} - e^{2t} \Delta u  + m^2 v=0,\,\,\,u(s,x)=0,\,\,\,u_{t}(s,x)=f(s, x).
\end{equation}
 In order to derive semilinear results for  \eqref{CauchyProbN2}, it is not sufficient to use the linear estimates from section \ref{generalestimates}, but
in addition one has to derive  $L^q-L^2$ estimates for the  parameter dependent Cauchy problem \eqref{linearassociate}.

The equation in \eqref{linearassociate} may be written as model \eqref{CauchyProb1} with $\alpha=1$ and $\delta= -4m^2<0$.
 Applying the change of variable
\[
v(\tau,  x)=u(t,  x), \qquad \tau+1=e^t, \qquad \tau_s+1=e^s
\]
the Cauchy problem \eqref{linearassociate} takes the form
\begin{align} \label{ChangVarMass}
v_{\tau\tau} -  \Delta v  + \frac{1}{1+\tau} v_{\tau} + \frac{m^2}{(1+\tau)^2} v=0,\,\,\,v(\tau_s, x)=0,\,\,\,v_{\tau}(\tau_s,  x)=(1+\tau_s)^{-1}f(s,x).
\end{align} By using the representation of solutions given by \cite{NPR}, one may derive the following:
\begin{prop}\label{propadd}
 If $m>0$ and $f(s,\cdot) \in \mathcal{D}_q$, $q\in [1,2)$, then
		the solution to \eqref{linearassociate} satisfies the following estimates:
\begin{equation}\label{estparameter}
\| \partial_tK_1(t,s, \cdot)\ast f(s, \cdot)\|_{L^2} + e^{t}\|\nabla_x K_1(t,s, \cdot)\ast f(s, \cdot)  \|_{L^2} \lesssim e^{-\frac{(s-t)}{2}} \| f(s, \cdot)\|_{L^2}
\end{equation}
\begin{equation}\label{estparameter1}
 \| K_1(t,s, \cdot)\ast f(s, \cdot) \|_{L^2} \lesssim e^{\frac{(s-t)}{2}} d(t,s) \| f(s, \cdot) \|_{L^2\cap L^q},
	\end{equation}
		for all $t\geq s\geq 0$, where
		\begin{align*}
		 d(t,s)= \left\{ \begin{array}{cr}
		1  &\mbox{\,\,\, for \,\,\,}  n>\frac{q}{2-q},\\
		(t-s)^{\frac{2-q}{2q}}   &\mbox{\,\,\, for \,\,\,} n=\frac{q}{2-q}.\end{array} \right.
		\end{align*}
\end{prop}
		For the consideration of the semilinear models we shall
use $q=1$.
 \begin{proof} (Theorem \ref{nonlinear})
We define
\[ X= \{ u\in \mathcal C([0,\infty),H^1(\mathbb{R}^n))\cap C^1([0, \infty), L^2(\mathbb{R}^n): \ \|u\|_X<\infty\}, \]
with the norm
\[\|u\|_X=\sup_{t\geq0}\left( e^{t/2}\big(  d(t)^{-1}\,\|u(t,\cdot)\|_{L^2}+ \|\nabla_x u(t,\cdot)\|_{L^2}+ e^{-t}\|u_t(t,\cdot)\|_{L^2}\big)\right).\]
For any~$u\in X$ we define
\[ Pu(t,x) :=  u^{lin}(t, x) + Nu(t,x), \]
where $u^{lin}(t, x):= K_0(t,0, x) \ast_{(x)} u_0(x) + K_1(t,0, x) \ast_{(x)} u_1(x)$ and
\[ Nu(t,x)= \int_0^tE_1(t,s, x)\ast |u(s,x)|^p ds. \]
Thanks to the derived estimates in  Example \ref{exexp} and Section \ref{generalestimates}, the linear part~$u^{lin}$ of the solution is in~$X$, and
\[ \|u^{lin}\|_X \lesssim \|(u_0,u_1)||_{\mathcal{D}_1(\mathbb{R}^n)}. \]
For any~$u, v\in X$, the nonlinear part of the solution satisfy
 \begin{eqnarray}
\label{eq:contractiondissipation}
\|Nu\|_X &\lesssim &\|u\|_X^p\\
\|Nu-Nv\|_{X}
     &\lesssim &\|u-v\|_{X(t)} \bigl(\|u\|_{X}^{p-1}+\|v\|_{X}^{p-1}\bigr) \label{eq:contractiondissipation2}
\end{eqnarray}
Hence the operator $P$ maps~$X$ into itself
and the existence of a
unique global solution $u$ follows by contraction \eqref{eq:contractiondissipation2} and continuation argument for small data.

By using the derived linear estimates \eqref{estparameter} and \eqref{estparameter1}, we prove \eqref{eq:contractiondissipation}, but we omit the proof of \eqref{eq:contractiondissipation2}, since it is analogous to the proof of \eqref{eq:contractiondissipation}. Indeed, for $\ell+j=0,1$ it holds
\[ \|\partial_t^{\ell}\nabla_x^j Nu(t,\cdot)\|_{L^2} \lesssim \int_0^{t} e^{\frac{(2\ell-1)t}{2}}d(t)^{(1-(j+\ell))} e^{\left(\frac{1}{2}-(j+\ell)\right)s}
\||u(s,\cdot)|^p\|_{L^2\cap L^1}ds.\]
Using
\[\||u(s,\cdot)|^p\|_{L^2\cap L^1}=\||u(s,\cdot)|^p\|_{L^2}+\||u(s,\cdot)|^p\|_{L^1}=\|u(s,\cdot)\|_{L^{2p}}^p+ \|u(s,\cdot)\|_{L^p}^p\]
and  applying Gagliardo-Nirenberg inequality, for all $p\geq 2$  and $ k=1,2$ we get
\[\|u(s,\cdot)\|_{L^{kp}}^p\lesssim
\|  u(s,\cdot)\|_{L^2}^{p(1-\theta(kp))} \| \nabla_x u(s,\cdot)\|_{L^2}^{p\theta(kp)}\lesssim  e^{-\frac{ps}{2}}d(s)^{p(1-\theta(kp))}\|u\|_X^p,\]
provided that $ p\leq \frac{n}{n-2}$ if $n\geq 3 $, i.e., $\theta(kp)\in [0,1]$. Therefore
\begin{eqnarray*}
\|\partial_t^{\ell}\nabla_x^j Nu(t,\cdot)\|_{L^2} &\lesssim e^{\frac{(2\ell-1)t}{2}}d(t)^{(1-(j+\ell))}
\int_0^{t} e^{\left(\frac{1}{2}-\frac{p}{2}\right)s} d(s)^{(1-\theta(kp))p} ds\|u\|_X^p\\
&\lesssim e^{\frac{(2\ell-1)t}{2}}d(t)^{(1-(j+\ell))}\|u\|_X^p,
\end{eqnarray*}
for all $p> 1$. This concludes the proof.
\end{proof}

%%%%%%%%%%%%%%%%%%%%%%%%%%%%%%%%%%
%%%%%%%%%%%%%%%%%%%%%%%%%%%%%%%%%%%
%%%%%%%%%%%%%%%%%%%%%%%%%%%%%%

\section{Concluding remarks and Open problems}\label{open}

\begin{enumerate}

\item The range of admissible $p$ in Theorem \ref{nonlinear}  came from the use of Gagliardo-Nirenberg inequality, so it's quite likely a technical restriction relate to  the choice of function spaces we take for the data and solutions.

	\item Through this paper we assume $a\notin L^1$.    In \cite{ERnew} the authors studied models for the damped wave models with integrable  in time speed of propagation. In a forthcoming  paper we shall study  the Cauchy problem  \eqref{wavegeneral}
	in the case $a\in L^1$,  which includes the Klein-Gordon equation in de Sitter spacetime, an important   model that appear in Mathematical Cosmology.
	
	\item
	Consider the semilinear problem for the scale invariant model
	\begin{equation*}
	\begin{cases}
	u_{tt} - (1+t)^{2\ell} \Delta u  + \frac{\widetilde{\mu}^2}{(1+t)^2} u
=|u|^p \\
(u(0,x),u_t(0,x)) = (u_0(x), u_1(x)),
\end{cases}
	\end{equation*}
	where  $\ell > 0$ and $\widetilde{\mu} = \mu (\ell+1)>0$. Applying the change of variable
	\[ v(\tau, x)=u(t,x), \qquad \tau+1=\frac{(1+t)^{\ell+1}}{\ell+1},\]
	this model can be transformed into the following Cauchy problem
\begin{equation*}
	\begin{cases}
			v_{\tau\tau} -  \Delta v  + \frac{\alpha}{1+\tau} v_{\tau} + \frac{\mu^2}{(1+\tau)^2} v=\beta^{-2}(1+\tau)^{-2\alpha}|v|^p \\
			 (v(\tau_0,x),v_{\tau}(\tau_0,x) )= (u_0(x),(1+\tau_0)^{-\ell} u_1(x)),
	\end{cases}
	\end{equation*}
	with $\alpha = \frac{\ell}{\ell + 1}$, $\beta=(\ell+1)^{\frac{\ell}{\ell+1}}$ and $\tau_0=-\frac{\ell}{\ell+1}$.
	By using data in $L^1\cap L^2$, in Section \ref{generalestimates} one may observe some  improvement in the decay rate for the energy solutions.
	This hints to the possibility to derive global (in time) existence of small data  energy solutions to this model. However, due to the fact that the dissipation is non-effective,   the approach  used in Theorem \ref{nonlinear} does not bring any sharp result for the critical exponent.
It is really a challenging problem to derive sharp global existence results to this model.
% (see \cite{NPR}).
\end{enumerate}

%%%%%%%%%%%%%%%%%%%%%%%%%%%%%%%%%%%%%%%%%%%%%%%%%%%%%%%%%%%%%

\end{document}